\documentclass[aap,MSNbibl,nameyear,dvips]{arximspdf}
\usepackage{algorithm}
\usepackage{algpseudocode}
\usepackage{graphicx}
\doi{10.1214/13-AAP969} %kopijuoti is PTS
\volume{24}
\issue{5}
\pubyear{2014}
\firstpage{1994}
\lastpage{2032}

\makeatletter
\newcommand{\rrvert}{\vert}
\newcommand{\llvert}{\vert}
\newcommand{\iint}{\int\!\!\!\int}
\def\coloneqq{:=}
\def\eqqcolon{=:}
\def\xrightarrow{\longrightarrow}
\newtheorem{theorem}{Theorem}[section]

\newtheorem{corollary}[theorem]{Corollary}

\newtheorem{lemma}[theorem]{Lemma}

\newproclaim{condition}[theorem]{Condition}
\newproclaim{definition}[theorem]{Definition}
\newproclaim{example}[theorem]{Example}

\newproclaim{remark}[theorem]{Remark}

\newcommand{\f}[2]{\frac{#1}{#2}}
\newcommand{\eqref}[1]{(\ref{#1})}
\renewcommand{\epsilon}{\varepsilon}
\makeatother

\begin{document}
\begin{frontmatter}

\title{Simulation of forward-reverse stochastic representations
for conditional diffusions}
\runtitle{Simulation of conditional diffusions}

\begin{aug}
\author{\fnms{Christian} \snm{Bayer}\corref{}\ead[label=e1]{christian.bayer@wias-berlin.de}}
\and
\author{\fnms{John} \snm{Schoenmakers}\ead[label=e2]{john.schoenmakers@wias-berlin.de}\thanksref{t2}}
\affiliation{WIAS Berlin}
\thankstext{t2}{Supported in part by the DFG Research Center \textsc{Matheon}
``Mathematics for Key Technologies'' in Berlin.}
\runauthor{C.~Bayer and J.~Schoenmakers}
\address{Weierstrass Institute\\
WIAS Berlin\\
Mohrenstr. 39\\
10117 Berlin\\
Germany\\
\printead{e1}\\
\phantom{E-mail:\ }\printead*{e2}}

\end{aug}

% HISTORY:
\received{\smonth{2} \syear{2013}}
\revised{\smonth{7} \syear{2013}}

% ABSTRACT
%
\begin{abstract}
In this paper we derive stochastic representations for the finite
dimensional distributions of a multidimensional diffusion on a fixed time
interval, conditioned on the terminal state. The conditioning can be with
respect to a fixed point or more generally with respect to some
subset. The
representations rely on a reverse process connected with the given (forward)
diffusion as introduced in Milstein, Schoenmakers and Spokoiny
[\textit{Bernoulli} \textbf{10} (2004) 281--312] in the context of a forward-reverse transition
density estimator. The
corresponding Monte Carlo estimators have essentially root-$N$ accuracy,
and hence
they do not suffer from the curse of dimensionality. We provide a detailed
convergence analysis and give a numerical example involving the realized
variance in a stochastic volatility asset model conditioned on a fixed
terminal value of the asset.
\end{abstract}

% KEYWORDS
%

%
% Pirmas kwd is didziosios raides

\end{frontmatter}

%s1 #&#
\section{Introduction}\label{sec1}

The central result in this paper is the development of a new generic procedure
for simulation of conditioned diffusions, also called diffusion bridges or
pinned diffusions. More specifically, for some given \mbox{(unconditional)} diffusion
process $X$, we aim to simulate the functional
%
%e1 #&#
\begin{equation}
\label{cp}\mathbb{E} \bigl[ g\bigl(X(s_{1}), \ldots,
X(s_{R})\bigr) \rrvert X(T)\in A, X(0) = x \bigr],
\end{equation}
where $0\leq s_{1}<s_{2}<\cdots<s_{R}<T$, $A$ is some set that may
consist of only one point, $g$ is an arbitrarily given suitable test
function and $x \in\mathbb{R}^d$ is a given state. In recent years, the
problem of computing terms such as \eqref{cp} has attracted a lot of attention
in the literature, sparked by several applications. Indeed, many relevant
properties of a diffusion process $X$ can be advantageously analyzed by
considering the process conditioned on certain appropriate events. One so
allows ``to study rare events by conditioning on the event happening or to
analyze the behaviour of a composite system when only some of its components
can be observed,'' as is eloquently put by \citet{HSV1}. For instance, in
statistical inference based on a continuous time model, discrete time
observations can be enriched to continuous time observations by
sampling from the
diffusion bridges between the discrete time data; see \citet{LCM} and
\citet{BlS} for more information. Conditional diffusions have further been
successfully used for critical calculations in rare event situations.
As an
example from computational chemistry, we refer to the review paper of
\citet{BCDG}, where diffusion bridges are used for detection of the transition
state surface between two stable regions $A$ and $B$ in configuration
space. Here, standard Monte Carlo simulation is prohibitively costly,
as the
event of such a transition is rare, provided that the ``walls'' in the energy
surface between $A$ and $B$ are high. However, by studying the process
conditioned on starting in $A$ and ending in $B$, one can efficiently observe
on which paths the configuration typically travels from $A$ to $B$. Other
possible applications appear in the field of stochastic environmental models,
for instance, regarding the concentration evolutions of pollution in
water; for
example see \citet{SHS} and references therein for a related problem.

Several approaches for simulation of diffusion bridges have already been
studied in the literature. For the theory of diffusion bridges we refer to
\citet{LZ} and the references therein. Many existing approaches utilize known
Radon--Nikodym densities of the law of the diffusion $X$ conditioned on initial
and terminal values, with respect to the law of a standard diffusion bridge
process (e.g., Wiener bridge) on path-space [as a Radon--Nikodym derivative
obtained by Doob's h-transform; see, e.g., \citet{RW} or
\citet{LZ}]. Several other approaches are based on (partial) knowledge
of the
transition densities of the unconditional diffusion (that is not generically
available, of course). For an overview of many different techniques, we refer
to \mbox{\citet{LCM}}.

First, let us mention the work by \citet{BPR} who construct a general,
rejection-based algorithm for solutions of \emph{one-dimensional} SDEs, based
on the Radon--Nikodym derivative of the law of the solution with
respect to the
Wiener measure. The algorithm gives (in finite, but random time) discrete
samples of the exact solution of the SDE. A~simple adaption of this algorithm
gives samples of the exact diffusion process conditioned on $X(T) = y$, by
using the law of the corresponding Brownian bridge as reference measure
(instead of the Wiener measure). An overview of related importance sampling
techniques is given by \citet{PR}. On the other hand, by relying on knowledge
of the transition densities of $X$, \citet{LCM} use a sequential
weighted Monte
Carlo framework, including resampling with optimal priority scores.

Another general technique used for simulation of diffusion bridges is the
Markov chain Monte Carlo method. Indeed, \citet{SVW} and \citet{HSV1}
show how
the law of a (multi-dimensional, uniformly elliptic, additive-noise) diffusion
$X$ conditioned on $X(T) = y$ can be regarded as the invariant
distribution of
a stochastic differential equation of Langevin type on path-space, that
is, of a
Langevin-type stochastic partial differential equation (SPDE). Thus, in
principle MCMC methods are applicable as explored by \citet{SVW}
and \citet{BRSV}. However, this requires the numerical solution of the SPDE
involved. It should be noted that in \citet{HSV2} the uniform ellipticity
condition is relaxed leading to a fourth order parabolic SPDE rather
than a
second order one.

Other notable approaches include those of \citet{MT}, which treat the
case of
physically relevant functionals of Wiener integrals with respect to Brownian
bridges, and \citet{S}, who uses an MCMC approach based on successive
modifications of the drift of the diffusion process.

Another approach is the one of \citet{BlS} developed for one-dimensional
diffusions. In order to obtain a sample from the process $X$
conditioned on
$X(0) = x$ and $X(T) = y$, \citet{BlS} start a path of the diffusion from
$(0,x)$ and another path of the diffusion in \emph{reversed time} at
$(T,y)$. If these paths hit at time $\tau$, consider the concatenated path
$Z$. The distribution of the process $Z$ (conditional on $0 \le\tau
\le T$)
equals the distribution of the bridge conditional on being hit by an
independent path of the underlying diffusion with initial distribution $p(0,
y, T, \cdot)$. As proved by \citet{BlS}, the probability of this event
approaches $1$ when $T \to\infty$. Finally, in order to improve the accuracy,
$Z$ is used as initial value of an MCMC algorithm on path space,
converging to
a sample from the true diffusion bridge.

A more general approach is given by \citet{DH} which relies on the explicit
Radon--Nikodym derivative of the diffusion $X$ conditioned on its
initial and
terminal values and another diffusion $Y$, which is modeled like the Brownian
bridge. In fact, $Y$ has the same dynamics as $X$, except for an extra term
$-\frac{Y(t) - y}{T-t}$ in the drift, which enforces $Y(T) = y$. Under certain
regularity conditions---in particular invertibility of the diffusion matrix
$\sigma= \sigma(t,x)$---\citet{DH} provide a Girsanov-type theorem, which
leads to a representation of the form
\[
\mathbb{E} \bigl[ g(X) | X(0) = x, X(T) = y \bigr] = \mathbb{E} \bigl[ g(Y) Z(Y)
\bigr]
\]
for functionals $g$ defined on path-space and a factor $Z(Y)$
explicitly given
as a functional of the path $Y$ together with quadratic variations of
functions of $Y$. As such this approach allows for direct Monte Carlo
simulation of \eqref{cp}. However, we stress that $Z(Y)$ explicitly
depends on
$\sigma^{-1}$ which does not exist in many hypo-elliptic applications.
On the
other hand, simulation of the bridge-type process $Y$ is numerically
troublesome because of the exploding drift term.

The new method presented in this article is inspired by the forward-reverse
estimator for the transition density $p(0,x,T,y)$ constructed by
\citet{MSS1}. Given a grid $0 \le s_0 < s_1 < \cdots< s_K = t^\ast<
t_1 <
\cdots< t_L = T$, we prove that
\[
\mathbb{E} \bigl[ g\bigl(X(s_1), \ldots, X(s_K),
X(t_1), \ldots, X(t_{L-1})\bigr) | X(s_0) =
x, X(T) = y \bigr]
\]
equals
%
%e2 #&#
\begin{equation}
\label{eq:for-rev-intro} \lim_{\epsilon\to0} \frac{\mathbb{E} [ g (X(s_1), \ldots, X(s_K),
Y(\widehat{t}_{L-1}), \ldots, Y(\widehat{t}_1) )
K_\epsilon ( Y(\widehat{t}_L) - X(t^\ast)  )
\mathcal{Y}(\widehat{t}_L)  ]}{
\mathbb{E} [ K_\epsilon ( Y(\widehat{t}_L) - X(t^\ast)  )
\mathcal{Y}(\widehat{t}_L)  ]},\hspace*{-20pt}
\end{equation}
which can be implemented by Monte Carlo simulation for any $\epsilon> 0$.
In \eqref{eq:for-rev-intro} $s_0 < t^\ast< T$ is a given grid-point
chosen by
the user. The process $X$ solves the original SDE with initial value
$X(s_0) =
x$ on the time-interval $[s_0, t^\ast]$. On the other hand, $Y$ is an
(independent) \emph{reverse process} as defined in Section~\ref{main} started
at $Y(0) = y$ and simulated until time $\widehat{t}_L \coloneqq T -
t^\ast$,
not on the original grid, but on a ``perturbed'' grid defined
in \eqref{eq:hat-grid}. [Note that $Y$ is different from the time-reversed
diffusion in the sense of \citet{HP} that explicitly requires the transition
density of $X$.] Indeed, the dynamics of $Y$ are explicitly given below in
terms of the dynamics of $X$---not relying on the transition density---and,
usually, share the same regularity properties; see \eqref{MC12}
and \eqref{revc}. Next, we weight the trajectories according to the distance
between $X(t^\ast)$ and $Y(T)$ using a kernel $K$ with bandwidth
$\epsilon$. Finally, we have an exponential weighting factor $\mathcal{Y}$,
similar to the Radon--Nikodym derivative in \citet{DH}. The denominator
in~\eqref{eq:for-rev-intro} actually corresponds to the forward-reverse
estimator for the transition density $p(s_0,x,T,y)$ of $X$ introduced by
\citet{MSS1}. The details of the Monte Carlo simulation are spelled out in
Section~\ref{analysis}, but we note that \eqref{eq:for-rev-intro} can be
computed to an accuracy of $\varepsilon$ with a complexity of
$\mathcal{O}(\varepsilon^{-2})$ in any dimension less or equal to
four\setcounter{footnote}{1}\footnote{In fact, this restriction can be lifted by use of higher order
kernels.} (disregarding possible discretization errors due to the
construction of samples $X$, $Y$ and~$\mathcal{Y}$).\footnote{The constant
will increase in the dimension. Moreover, we ignore the cost of checking
equality of two integers.} Thus our algorithm essentially achieves the
optimal rate of convergence for Monte Carlo algorithms.

We underline that the forward-reverse algorithm for \eqref{cp}
presented here
is not a straightforward extension of the forward-reverse algorithm for
transition densities of \citet{MSS1}. The main difficulty lies in the extension
of the representation from just one intermediate time $0 < t^\ast< T$
to an
arbitrary time grid $0 \le s_0 < s_1 < \cdots< s_K = t^\ast< t_1 <
\cdots<
t_L = T$ with $L>1$. In the nonautonomous case this issue is further
complicated due to the fact that the dynamics of the reverse process as
defined in \citet{MSS1} depends explicitly on both $t^\ast$ and
$T$. Obviously, the
different structure also requires a different error analysis. In
particular, we
need sharper error bounds than \citet{MSS1}.

In comparison to the other methods mentioned above, our new procedure
has the
following main features:

\begin{longlist}[(iii)]
\item[(i)] The method applies to multidimensional diffusions.

\item[(ii)] It is based on simulation of \emph{unconditional}
diffusions only,
hence technical simulation problems due to exploding drifts in SDEs that
govern particular diffusion bridges are avoided.

\item[(iii)] The vector fields determining the (forward) SDE that
governs $X$
only need to satisfy a H\"ormander-type condition guaranteeing sufficient
regularity and exponential decay of the transition densities. In
particular, the
diffusion matrix of $X$ may be degenerate.

\item[(iv)] The estimator corresponding to the developed stochastic
representation for~(\ref{cp}) is root-$N$ consistent, that is the mean
square estimation accuracy is of order $\mathcal{O}(N^{-1/2})$ with $N$ being
the number of trajectories that need to be simulated.
\end{longlist}

As a matter of fact, the methods for simulating diffusion bridges known
in the
literature so far, do not cover all the features (i)--(iv)
simultaneously. For
example, \citet{DH}
require that either the diffusion matrix is invertible, or impose some very
specific structural conditions on the drift and diffusion matrix of the
process $X$. Moreover, the exploding drift terms in their process $Y$ makes
simulation of the auxiliary process $Y$ nontrivial. On the other hand, the
method of \citet{BlS} in germ carries some ideas related to our
approach, but
they need to impose balance restrictions on the transition density of
$X$, and
moreover their method---together with several others---is only
one dimensional. The methods of \citet{SVW} and the related papers mentioned
above also involve some further structural assumptions and, in addition,
require numerical solutions of SPDEs.

Moreover, we complement our algorithm by an adaptation, which allows us to
treat the more general problem of conditioning at final time $T$ not on all,
but just on some components of the vector $X_T$. More precisely, we
present a
variant of the algorithm for computing conditional expectations where
$X_T$ is
conditioned to lie in a ``simple'' set $A$, that is, either $A$ has positive
measure both under the Lebesgue measure and the distribution of $X_T$,
or $A$
is an affine plane of dimension $0 \leq d' \leq d$. In order to achieve this
extension, we need to prove (Lebesgue) integrable error bounds for the
forward-reverse algorithm for the case where $X_T$ is conditioned to a value~$y$.

The structure of the paper is as follows. In Section~\ref{recap} we
recap the
essential facts concerning the reverse diffusion system of \citet
{MSS1}. The
main representation theorems for the diffusion conditioned on reaching
a fixed
state, or conditioned on reaching some Borel set, are derived in
Section~\ref{main}. A detailed accuracy analysis concerning the Monte Carlo
estimators for the respective conditioned diffusions is provided in
Section~\ref{analysis}, including the precise required regularity assumptions
given in Conditions \ref{ass:bound-density}, \ref{ass:kernel-order}
and \ref{ass:convenience}. Limitations of the method are discussed in
Section~\ref{sec:limit-forw-reverse}, while Section~\ref{num} provides a
numerical study involving a Heston-type stochastic volatility model.

%s2 #&#
\section{Recap of forward-reverse representations for diffusions}
\label{recap}

In this section we recapitulate shortly the main ingredients in the approach
by \citet{MSS1}. Let us consider the SDE
%
%e3 #&#
\begin{eqnarray}\label{In1}
dX_{t,x}(s)=a\bigl(s,X_{t,x}(s)\bigr)\,ds+\sigma
\bigl(s,X_{t,x}(s)\bigr)\,dW(s),
\nonumber
\\[-8pt]
\\[-8pt]
 \eqntext{0\leq s\leq T, X_{t,x}(t)=x,}
\end{eqnarray}
where $X_{t,x}\in\mathbb{R}^{d}$, $a\dvtx [t,T]\times\mathbb
{R}^{d}\rightarrow
\mathbb{R}^{d}$, $\sigma\dvtx [t,T]\times\mathbb{R}^{d}\rightarrow\mathbb
{R}%
^{d\times m}$, $W$ is an $m$-dimensional standard Wiener process and
$x\in
\mathbb{R}^{d}$. At this stage, we only assume that $X$ admits a $C^{2}$
transition density $p$ and that the coefficients of \eqref{In1} are $C^{2}$
as well.

Along with the (forward) process $X$ given by (\ref{In1}), \citet{MSS1}
introduced an associated process $(Y_{t,y;T}(s), \mathcal
{Y}_{t,y;T}(s))$ in $\mathbb{R}^{d}\times\mathbb{R}$, $t\leq s\leq
T$, termed \emph{reverse} process on the interval $[t,T]$, that solves the
SDE
%
%e4 #&#
\begin{eqnarray}\label{MC12}
dY_{t,y;T}(s)&=&\alpha_{t,T}
\bigl(s,Y_{t,y;T}(s)\bigr) \,ds \nonumber\\
&&{}+ \widetilde{\sigma}_{t,T}
\bigl(s,Y_{t,y;T}(s)\bigr) \,d\widetilde{W}(s),\qquad Y_{t,y;T}(t)=y,
\\
\mathcal{Y}_{t,y;T}(s)&=&\exp \biggl( \int_{t}^{s}c_{t,T}
\bigl(u,Y_{t,y;T}(u)\bigr)\,du \biggr) \nonumber
\end{eqnarray}
with $\widetilde{W}$ being a (from $W$ independent) $m$-dimensional Wiener
process, and
%
%e5 #&#
\begin{eqnarray} \label{revc}
\alpha_{t,T}^{i}(s,y)& \coloneqq&\sum
_{j=1}^{d}\frac{\partial}{\partial
y^{j}}b^{ij}(T+t-s,y)-a^{i}(T+t-s,y),\qquad
b\coloneqq\sigma\sigma ^{\top},
\nonumber
\\
\widetilde{\sigma}_{t,T}(s,y)& \coloneqq&\sigma(T+t-s,y),
\\
\qquad c_{t,T}(s,y)& \coloneqq&\frac{1}{2}\sum
_{i,j=1}^{d}\frac{\partial
^{2}b^{ij}}{%
\partial y^{i}\partial y^{j}}(T+t-s,y)-\sum
_{i=1}^{d}\frac{\partial
a^{i}}{%
\partial y^{i}}(T+t-s,y).
\nonumber
\end{eqnarray}
Despite its name, we stress that $(Y,\mathcal{Y})$ is the solution of an
ordinary SDE \emph{forward} in time on the interval $[t,T]$.

One of the central results in \citet{MSS1} is the following theorem.

%th2.1 #&#
\begin{theorem}[{[M.S.S. (2004)]}]
\label{MT} For fixed $t,x,y$ and $t<t^{\ast}<T$, and any bi-variate test
function $f  $we have
%
%e6 #&#
\begin{eqnarray}\label{Jf}
&&\mathbb{E} \bigl[f\bigl(X_{t,x}\bigl(t^{\ast}
\bigr),Y_{t^{\ast},y;T}(T)\bigr) \mathcal{Y}%
_{t^{\ast},y;T}(T)\bigr]
\nonumber
\\[-8pt]
\\[-8pt]
\nonumber
&&\qquad=
\iint p\bigl(t,x,t^{\ast},x^{\prime}\bigr)p\bigl(t^{\ast
},y^{\prime
},T,y
\bigr)f\bigl(x^{\prime},y^{\prime}\bigr)\,dx^{\prime}\,dy^{\prime},
\end{eqnarray}
where $ X_{t,x}(s) $ satisfies the forward equation (\ref{In1}), and $%
(Y_{t^{\ast},y;T}(s),\mathcal{Y}_{t^{\ast},y;T}(s))$, $s\geq t^{\ast}$,
is the solution of the reverse system (\ref{MC12}).
\end{theorem}

%co2.2 #&#
\begin{corollary}
\label{cor:MT} By taking $f\equiv1$, (\ref{Jf}) yields
%
%e7 #&#
\begin{equation}
\mathbb{E} \bigl[\mathcal{Y}_{t^{\ast},y;T}(T)\bigr]=\int p
\bigl(t^{\ast},y^{\prime
},T,y\bigr)\,dy^{\prime}, \label{Jf1}
\end{equation}
which obviously extends to $t^{\ast}=t$. By next taking $f(x^{\prime
},y^{\prime})=f(x^{\prime})$ (while abusing notation slightly) we obtain
from (\ref{Jf}), using (\ref{Jf1}) and the independence of $X$ and $(Y,%
\mathcal{Y})$,
\[
\mathbb{E} \bigl[f\bigl(X_{t,x}\bigl(t^{\ast}\bigr)\bigr)
\bigr]=\int p\bigl(t,x,t^{\ast},x^{\prime
}\bigr)f
\bigl(x^{\prime}\bigr)\,dx^{\prime},
\]
which obviously extends to $t^{\ast}=T$, that is, the standard forward
stochastic representation for $\int p(t,x,T,x^{\prime})f(x^{\prime
})\,dx^{\prime}$. On the other hand, by taking $f(x^{\prime},y^{\prime
})=f(y^{\prime})$ we obtain the so called reverse stochastic
representation%
%
%e8 #&#
\begin{equation}
\mathbb{E} \bigl[f\bigl(Y_{t^{\ast},y;T}(T)\bigr) \mathcal{Y}_{t^{\ast},y;T}(T)
\bigr]=\int p\bigl(t^{\ast},y^{\prime},T,y\bigr)f
\bigl(y^{\prime}\bigr)\,dy^{\prime}, \label{MT1}
\end{equation}
which obviously extends to $t^{\ast}=t$.
\end{corollary}

%s3 #&#
\section{Forward-reverse representations for conditional diffusions}
\label{main}

It should be noted that in \citet{MSS1} the time domain of the reverse
process was considered fixed. For our purposes however, it turns out to be
more effective (in particular regarding the proof of Theorem \ref{key}
below) to consider reverse processes suitably defined on different time
domains. In particular it turns out be fruitful to formulate the
forward-reverse representations of the previous section in terms of reverse
processes defined on $[0,T]$ for suitable $T>0$. We therefore introduce
the\vadjust{\goodbreak}
reverse\break  process
%
%e9 #&#
\begin{equation}
\bigl(Y_{y;T}(s),\mathcal{Y}_{y;T}(s)\bigr)_{0\leq s\leq T}
\coloneqq \bigl(Y_{0,y;T}(s),\mathcal{Y}_{0,y;T}(s)
\bigr)_{0\leq s\leq T} \label{so}
\end{equation}
that starts at time $s=0$ at a generic state $ ( y,1 ) $, is
defined on an interval $[0,T]$ and satisfies (\ref{MC12}) with
coefficients (\ref{revc}) for $t=0$, that is, (\ref{so}) solves the SDE%
%
%e10 #&#
\begin{eqnarray}\label{so1}
dY(s)&=&\alpha_{0,T}\bigl(s,Y(s)\bigr)\,ds+
\widetilde{\sigma }_{0,T}\bigl(s,Y(s)\bigr)\,d\widetilde{W}%
(s), \qquad Y(0)=y,
\nonumber
\\[-8pt]
\\[-8pt]
\nonumber
\mathcal{Y}(s)&=&\exp \biggl( \int_{0}^{s}c_{0,T}
\bigl(u,Y(u)\bigr)\,du \biggr).
\end{eqnarray}
As a result, we have for any fixed $t$, $0\leq t\leq T$, that
\[
\bigl( Y_{y;T}(s),\mathcal{Y}_{y;T}(s) \bigr)
_{0\leq s\leq T-t}= \bigl( Y_{t,y;T}(t+s),\mathcal{Y}_{t,y;T}(t+s)
\bigr) _{0\leq s\leq T-t},
\]
whence (\ref{Jf}) and (\ref{MT1}) may be equivalently written as%
%
%e11 #&#
%e12 #&#
\begin{eqnarray}
\label{Jf1a} &&\mathbb{E} \bigl[f\bigl(X_{t,x}\bigl(t^{\ast}
\bigr),Y_{y;T}\bigl(T-t^{\ast}\bigr)\bigr)\mathcal{Y}%
_{y;T}\bigl(T-t^{\ast}\bigr) \bigr]
\nonumber
\\[-8pt]
\\[-8pt]
\nonumber
&&\qquad=
\iint p\bigl(t,x,t^{\ast},x^{\prime}\bigr)p\bigl(t^{\ast
},y^{\prime},T,y
\bigr)f\bigl(x^{\prime},y^{\prime}\bigr)\,dx^{\prime}\,dy^{\prime}
\end{eqnarray}
and%
%
%e13 #&#
\begin{equation}
\mathbb{E} \bigl[f\bigl(Y_{y;T}\bigl(T-t^{\ast}\bigr)\bigr)
\mathcal{Y}_{y;T}\bigl(T-t^{\ast}\bigr) \bigr]=\int p
\bigl(t^{\ast},y^{\prime},T,y\bigr)f\bigl(y^{\prime}
\bigr)\,dy^{\prime}, \label{MT1a}
\end{equation}
respectively. The main benefit is that the reverse process used in
representations~(\ref{Jf}) and (\ref{MT1}) depend on both $t^{\ast}$
and $T$,\vadjust{\goodbreak} while the one used in (\ref{Jf1a}) and (\ref{MT1a}) depends on $T$ only.
In particular, (\ref{MT1a}) may be considered as a reverse representation
for all $0<t^{\ast}<T$ in terms of one and the same reverse process $%
(Y_{y;T},\mathcal{Y}_{y;T})$.

%s3.1 #&#
\subsection{Representations for conditioning on a fixed state}
\label{sec:repr-cond-fixed}

Let us start with the following lemma.

%le3.1 #&#
\begin{lemma}
\label{lem} For any $0<s<t\leq T$ it holds that
\begin{eqnarray*}
Y_{Y_{y;T}(s);T-s}(t-s) &=&Y_{y;T}(t),
\\
\mathcal{Y}_{y;T}(t) &=&\mathcal{Y}_{y;T}(s)
\mathcal{Y}%
_{Y_{y;T}(s);T-s}(t-s).
\end{eqnarray*}
\end{lemma}

\begin{pf}
The first statement is directly obvious from (\ref{so1}). From this the
second statement follows by
\begin{eqnarray*}
\mathcal{Y}_{y;T}(t) &=&e^{\int_{0}^{s}c_{0,T}(u,Y_{y;T}(u))\,du}e^{%
\int_{s}^{t}c_{0,T}(u,Y_{Y_{y;T}(s);T-s}(u-s))\,du}
\\
&=&\mathcal{Y}_{y;T}(s)e^{\int_{0}^{t-s}c_{0,T-s}(u,Y_{Y_{y;T}(s);T-s}(u))\,du}
\\
&=&\mathcal{Y}_{y;T}(s)\mathcal{Y}_{Y_{y;T}(s);T-s}(t-s).%\qedhere
\end{eqnarray*}
\upqed\end{pf}\eject

We are now ready to state the following key theorem.

%th3.2 #&#
\begin{theorem}
\label{key} Given a grid $\mathcal{D}_{L}:=\{0\leq t^{\ast
}<t_{1}<\cdots
<t_{L}\}$, it holds that
\begin{eqnarray*}
&&\mathbb{E} \bigl[ f\bigl(Y_{y;t_{L}}(t_{L}-t_{0}),Y_{y;t_{L}}(t_{L}-t_{1}),
\ldots,Y_{y;t_{L}}(t_{L}-t_{L-1})\bigr)
\mathcal{Y}_{y;t_{L}}(t_{L}-t_{0}) \bigr]
\\
&&\qquad= \int_{\mathbb{R}^{d\times
L}}f(y_{0},y_{1},
\ldots,y_{L-1})\prod_{i=1}^{L}p(t_{i-1},y_{i-1},t_{i},y_{i})\,dy_{i-1}
\end{eqnarray*}
with $y_{L} \coloneqq y$ and $t_{0} \coloneqq t^{\ast}$.
\end{theorem}

\begin{pf}
The following proof---much shorter than the proof in a previous
version of
the paper---has essentially been pointed out to us by an anonymous referee.
We fix $t_{0}$ ($=t^{\ast}$) and use induction on $L$. For $L=1$ the
statement boils down to~(\ref{MT1a}) with $T:=t_{1}$. Suppose the statement
is proved for some $L\geq1$. For the grid
\[
\mathcal{D}_{L+1}=\{t_{0}<t_{1}<
\cdots<t_{L+1}\}
\]
we next consider for any test function $f\dvtx \mathbb{R}^{d\times
(L+1)}\rightarrow\mathbb{R}$,
\begin{eqnarray*}
&&\mathbb{E} \bigl[ f\bigl(Y_{y;t_{L+1}}(t_{L+1}-t_{0}),Y_{y;t_{L+1}}(t_{L+1}-t_{1}),
\ldots,Y_{y;t_{L+1}}(t_{L+1}-t_{L})\bigr)\\
&&\hspace*{202pt}{}\times \mathcal
{Y}_{y;t_{L+1}}(t_{L+1}-t_{0}) \bigr]
\\
&&\qquad=\mathbb{E} \bigl[ \mathbb{E} \bigl[ f \bigl( Y_{Y_{y;t_{L+1}}(t_{L+1}-t_{L});t_{L}}(t_{L}-t_{0}),
Y_{Y_{y;t_{L+1}}(t_{L+1}-t_{L});t_{L}}(t_{L}-t_{1}), \ldots,
\\
&&\hspace*{68pt}\qquad\quad Y_{Y_{y;t_{L+1}}(t_{L+1}-t_{L});t_{L}}(t_{L}-t_{L-1}),
Y_{Y_{y;t_{L+1}}(t_{L+1}-t_{L});t_{L}}(0) \bigr)\\
&&\hspace*{24pt}\qquad\quad{} \times
\mathcal {Y}_{y;t_{L+1}}(t_{L+1}-t_{L})
\mathcal{Y}_{Y_{y;t_{L+1}}(t_{L+1}-t_{L});t_{L}}(t_{L}-t_{0}) \vert\\
&&\hspace*{166pt}
Y_{y;t_{L+1}}(t_{L+1}-t_{L}), \mathcal{Y}_{y;t_{L+1}}(t_{L+1}-t_{L})
\bigr] \bigr]
\\
&&\qquad=\mathbb{E} \bigl[ \mathcal{Y}_{y;t_{L+1}}(t_{L+1}-t_{L})
\\
&&\hspace*{43pt}{}\times\mathbb {E} \bigl[ f \bigl(Y_{Y_{y;t_{L+1}}(t_{L+1}-t_{L});t_{L}}(t_{L}-t_{0}),
Y_{Y_{y;t_{L+1}}(t_{L+1}-t_{L});t_{L}}(t_{L}-t_{1}), \ldots,
\\
&&\hspace*{93pt} \qquad\quad Y_{Y_{y;t_{L+1}}(t_{L+1}-t_{L});t_{L}}(t_{L}-t_{L-1}),
Y_{y;t_{L+1}}(t_{L+1}-t_{L})\bigr)
\\
& &\hspace*{97pt}\qquad\quad{}\times\mathcal{Y}_{Y_{y;t_{L+1}}(t_{L+1}-t_{L});t_{L}}(t_{L}-t_{0}) \vert
Y_{y;t_{L+1}}(t_{L+1}-t_{L}) \bigr] \bigr].
\end{eqnarray*}
By the induction hypothesis, we have that
\begin{eqnarray*}
&&\mathbb{E} \bigl[ f\bigl(Y_{Y_{y;t_{L+1}}(t_{L+1}-t_{L});
t_{L}}(t_{L}-t_{0}),Y_{Y_{y;t_{L+1}}(t_{L+1}-t_{L});t_{L}}(t_{L}-t_{1}),
\ldots,\\
&&\hspace*{70pt}Y_{Y_{y;t_{L+1}}(t_{L+1}-t_{L});t_{L}}(t_{L}-t_{L-1}),
Y_{y;t_{L+1}}(t_{L+1}-t_{L})\bigr)\\
&&\hspace*{55pt}{}\times \mathcal
{Y}_{Y_{y;t_{L+1}}(t_{L+1}-t_{L});t_{L}}(t_{L}-t_{0}) \vert
Y_{y;t_{L+1}}(t_{L+1}-t_{L})=z \bigr]
\\
&&\qquad= \mathbb{E} \bigl[ f\bigl(Y_{z;t_{L}}(t_{L}-t_{0}),
Y_{z;t_{L}}(t_{L}-t_{1}),
\ldots,Y_{z;t_{L}}(t_{L}-t_{L-1}),z\bigr)
\mathcal{Y}_{z;t_{L}}(t_{L}-t_{0})%
\bigr]
\\
&&\qquad= \int_{\mathbb{R}^{d\times L}}f(y_{0},y_{1},
\ldots,y_{L-1},z) \prod_{i=1}^{L}p(t_{i-1},y_{i-1},t_{i},y_{i})
\,dy_{i-1} \eqqcolon F(z)
\end{eqnarray*}
with $y_L \coloneqq z$, and so we obtain
\begin{eqnarray*}
&&\mathbb{E} \bigl[ f\bigl(Y_{y;t_{L+1}}(t_{L+1}-t_{0}),Y_{y;t_{L+1}}(t_{L+1}-t_{1}),
\ldots,Y_{y;t_{L+1}}(t_{L+1}-t_{L})\bigr)\\
&&\hspace*{202pt}{}\times \mathcal
{Y}_{y;t_{L+1}}(t_{L+1}-t_{0}) \bigr]
\\
&&\qquad=\mathbb{E} \bigl[ \mathcal{Y}_{y;t_{L+1}}(t_{L+1}-t_{L})F
\bigl( Y_{y;t_{L+1}}(t_{L+1}-t_{L}) \bigr) \bigr]
\\
&&\hspace*{-3pt}\qquad\stackrel{\scriptsize{(\ref{MT1a})}}{=} \int p(t_{L},z,t_{L+1},y)F(z)\,dz
\\
&&\qquad=\int_{\mathbb{R}^{d\times ( L+1 ) }}f(y_{0},y_{1},\ldots
,y_{L})\prod_{i=1}^{L+1}p(t_{i-1},y_{i-1},t_{i},y_{i})\,dy_{i-1},
\end{eqnarray*}
where $y_{L+1} \coloneqq y$ and the integration variable $z$ is renamed to
$y_{L}$.
\end{pf}

For the next theorem, we consider an extended time grid
%
%e14 #&#
\begin{equation}
\label{eq:full-grid} 0 \le s_0 < s_1 < \cdots<
s_K = t^\ast= t_0 < t_1 <
\cdots< t_L = T.
\end{equation}
For convenience, we also introduce the notation
%
%e15 #&#
\begin{equation}
\label{eq:hat-grid} \widehat{t}_i \coloneqq t_L -
t_{L-i},\qquad  i=1, \ldots, L.
\end{equation}
Moreover, we fix a starting point $x\in\mathbb{R}^{d}$.

%th3.3 #&#
\begin{theorem}
\label{Gen1} For any $f\dvtx \mathbb{R}^{d\times(K+L)}\rightarrow\mathbb{R}$ and
grids \eqref{eq:full-grid} together with \eqref{eq:hat-grid}, we have
\begin{eqnarray*}
&&\mathbb{E} \bigl[f\bigl(X_{s_{0},x}(s_{1}), \ldots,
X_{s_{0},x}(s_{K}), Y_{y;T}(\widehat{t}_L),
Y_{y;T}(\widehat{t}_{L-1}), \ldots, Y_{y;T}(
\widehat{t}_{1})\bigr) \mathcal{Y}_{y;T}(
\widehat{t}_L) \bigr]
\\
&&\qquad= \int_{\mathbb{R}^{d\times(K+L)}} f(x_{1},\ldots,x_{K},y_{0},y_{1},
\ldots,y_{L-1}) \\
&&\hspace*{75pt}{}\times\prod_{i=1}^{K}
p(s_{i-1},x_{i-1},s_{i},x_{i})\,dx_{i}
\prod_{i=1}^{L} p(t_{i-1},y_{i-1},t_{i},y_{i})\,dy_{i-1}
\end{eqnarray*}
with $x_{0}\coloneqq x$, $y_{L}\coloneqq y$, and the processes $X$ and
$(Y,\mathcal{Y})$ being independent.
\end{theorem}

Theorem \ref{Gen1} follows directly from Theorem \ref{key} by a standard
conditioning argument and the Chapman--Kolmogorov equation. Note that for
$K=L=1$, Theorem \ref{Gen1} collapses to Theorem \ref{MT}.

We are now ready to derive a forward-reverse stochastic representation
for the
finite dimensional distributions of the process $X_{s_{0},x}$,
conditional on
$X_{s_{0},x}(T)=y$, for fixed $s_{0}<T$, and fixed
$x,y\in\mathbb{R}^{d}$. To this end we henceforth assume that
%
%e16 #&#
\begin{equation}
p(s_0,x,T,y)>0.\label{assp}
\end{equation}
We also need to assume continuity of $p$. Let us take a bounded measurable
test function
\[
g(x_{1},\ldots,x_{K},y_{1},\ldots,y_{L-1})
\dvtx \mathbb{R}^{d\times
(K+L-1)}\rightarrow \mathbb{R},
\]
and consider the conditional expectation
%
%e17 #&#
\begin{eqnarray}\label{CE}
&&\mathbb{E} \bigl[  g\bigl(X_{s_{0},x}(s_{1}),\ldots,X_{s_{0},x}
(s_{K-1}),X_{s_{0},x}\bigl(t^{\ast}\bigr),X_{s_{0},x}(t_{1}),\ldots,X_{s_{0},x}
(t_{L-1})\bigr)\rrvert
\nonumber
\\[-8pt]
\\[-8pt]
\nonumber
&&\hspace*{257pt}{} X_{s_{0},x}(T)=y \bigr].
\end{eqnarray}
The distribution of the diffusion $X_{s_{0},x}$ conditional on $X_{s_{0},x}(T)
= y$ is completely determined by the totality of conditional
expectations of
the form (\ref{CE}). These conditional expectations may be obtained due to
Theorem \ref{thr:conditional-dist-general-grid} below.

%th3.4 #&#
\begin{theorem}
\label{thr:conditional-dist-general-grid} Consider the forward process $X$
and its reverse process $(Y, \mathcal{Y})$ as before and the grids as
specified in \eqref{eq:full-grid} and \eqref{eq:hat-grid}. Let
\[
K_{\epsilon}(u)\coloneqq\epsilon^{-d}K(u/\epsilon), \qquad y\in
\mathbb{R}^{d},
\]
with $K$ being integrable on $\mathbb{R}^{d}$ and $\int_{\mathbb{R}^{d}%
}K(u)\,du=1$. Hence, formally $K_{\epsilon}$ converges to the delta function
$\delta_{0}$ on $\mathbb{R}^{d}$ (in distribution sense) as $\epsilon
\downarrow0$. Then, since $p(s_0,x,T,y) > 0$ by assumption,\vadjust{\goodbreak} for any bounded
measurable function $g\dvtx \mathbb{R}^{d\times(K+L-1)} \to\mathbb{R}$,
we have
%
%e18 #&#
%e19 #&#
%e20 #&#
\begin{eqnarray}
\label{cdy1}&& \mathbb{E} \bigl[  g \bigl( X_{s_{0},x}(s_{1}),
\ldots, X_{s_{0},x}\bigl(t^{\ast}\bigr), X_{s_{0},x}(t_{1}),
\ldots, X_{s_{0},x}(t_{L-1}) \bigr) \rrvert X_{s_{0},x}(T)
= y \bigr] \nonumber\\
&&\qquad=
\frac{1}{p(s_{0},x,T,y)}
\nonumber
\\[-8pt]
\\[-8pt]
\nonumber
&&\qquad\quad{}\times\lim_{\epsilon\downarrow0} \mathbb{E} \bigl[ g
\bigl(X_{s_{0},x}(s_{1}), \ldots, X_{s_{0},x}
\bigl(t^{\ast}\bigr), Y_{y;T}(\widehat{t}_{L-1}),
\ldots, Y_{y;T}(\widehat{t}_{1})\bigr)
\\
&&\hspace*{113pt}\qquad\quad\times K_{\epsilon} \bigl( Y_{y;T}(\widehat{t}_L)-X_{s_{0},x}
\bigl(t^{\ast
}\bigr) \bigr) \mathcal{Y}_{y;T}(
\widehat{t}_L) \bigr].\nonumber
\end{eqnarray}
\end{theorem}

\begin{pf}
By applying Theorem \ref{Gen1} to
\[
f(x_{1}, \ldots, x_{K}, y_{0},
y_{1}, \ldots, y_{L-1}) \coloneqq g(x_{1},
\ldots, x_{K},y_{1}, \ldots, y_{L-1})
K_{\epsilon}(y_{0}-x_{K}),
\]
we obtain
%
%e21 #&#
\begin{eqnarray}\label{rh}
&&\mathbb{E} \bigl[g \bigl( X_{s_{0},x}(s_{1}), \ldots,
X_{s_{0},x}\bigl(t^{\ast}\bigr), Y_{y;T}(
\widehat{t}_{L-1}), \ldots, Y_{y;T}(\widehat{t}_{1})
\bigr) \nonumber\\
&&\hspace*{79pt}{}\times K_{\epsilon} \bigl( Y_{y;T}(\widehat{t}_L) -
X_{s_{0},x}\bigl(t^{\ast}\bigr) \bigr) \mathcal{Y}_{y;T}(
\widehat{t}_L) \bigr]
\nonumber
\\
&&\qquad= \int_{\mathbb{R}^{d\times(K+L)}} g(x_{1}, \ldots, x_{K},
y_{1},\ldots, y_{L-1}) K_{\epsilon}(y_{1}-x_{K})
\nonumber
\\[-8pt]
\\[-8pt]
\nonumber
&&\hspace*{42pt}\qquad\quad{} \times\prod_{i=1}^{K}
p(s_{i-1},x_{i-1},s_{i},x_{i})
\,dx_{i} \prod_{i=1}^{L}
p(t_{i-1},y_{i-1},t_{i},y_{i})
\,dy_{i-1}
\\
&&\qquad= \int_{\mathbb{R}^{d\times(K+L)}} g(x_{1}, \ldots,
x_{K},y_{1},\ldots, y_{L-1}) K(v)\,dv\, p
\bigl(t^\ast,x_k + \epsilon v,t_1,y_{2}
\bigr) \nonumber
\\
&&\hspace*{42pt}\qquad\quad{} \times\prod_{i=1}^{K}
p(s_{i-1}, x_{i-1}, s_{i}, x_{i})
\,dx_{i} \prod_{i=2}^{L}
p(t_{i-1},y_{i-1},t_{l},y_{i})\,dy_{i-1}.
\nonumber
\end{eqnarray}
By sending $\epsilon$ to zero, (\ref{rh}) clearly converges to
\begin{eqnarray*}
&&\int_{\mathbb{R}^{d\times(K+L-1)}}g(x_{1},\ldots,x_{K},y_{1},\ldots,y_{L-1})
\prod_{i=1}^{K} p(s_{i-1},x_{i-1},s_{i},x_{i})
\,dx_{i}
\\
&&\hspace*{28pt}\qquad{}\times p\bigl(t^{\ast},x_{K},t_{1},
y_{1}\bigr) \prod_{i=2}^{L}
p(t_{i-1},y_{i-1},t_{l},y_{i})\,dy_{i-1},
\end{eqnarray*}
from which (\ref{cdy1}) easily follows.
\end{pf}

If the original grid $t^{\ast}=t_{0}< \cdots<t_{L}=T$ is
equidistant, then the transformed grid $0 = \widehat{t}_0 < \cdots<
\widehat{t}_L = T-t^\ast$ is obtained by a translation with $-t^\ast$, which
leads to the following corollary.
%
%co3.5 #&#
\begin{corollary}
\label{cor-equidist-grid} If the time grid $t^{\ast}=t_{0}< \cdots<t_{L}=T$
is equidistant, we have
\begin{eqnarray*}
&&\mathbb{E} \bigl[  g \bigl( X_{s_{0},x}(s_{1}),
\ldots,X_{s_{0},x}\bigl(t^{\ast}\bigr),X_{s_{0},x}(t_{1}),
\ldots,X_{s_{0},x}(t_{L-1}) \bigr) \rrvert X_{s_{0},x}(T)=y
\bigr]
\\
&&\qquad=\frac{1}{p(s_{0},x,T,y)}\\
&&\qquad\quad{}\times\lim_{\epsilon\downarrow0}\mathbb{E} \bigl[ g
\bigl(X_{s_{0},x}(s_{1}), \ldots, X_{s_{0},x}
\bigl(t^{\ast
}\bigr),Y_{y;T}\bigl(t_{L-1}-t^\ast
\bigr), \ldots, Y_{y;T}\bigl(t_{1}-t^\ast\bigr)
\bigr)
\\
&&\hspace*{110pt}\qquad\quad{}\times K_{\epsilon} \bigl( Y_{y;T}\bigl(T-t^\ast
\bigr)-X_{s_{0},x}\bigl(t^{\ast}\bigr) \bigr) \mathcal{Y}_{y;T}
\bigl(T-t^\ast\bigr) \bigr].
\end{eqnarray*}
Moreover, by setting $g\equiv1$, we retrieve the forward-reverse
representation of the transition density in \citet{MSS1},
%
%e22 #&#
\begin{equation}
p(s_{0},x,T,y)=\lim_{\epsilon\downarrow0}\mathbb{E}
\bigl[K_{\epsilon
} \bigl( Y_{y;T}\bigl(T-t^\ast
\bigr)-X_{s_{0},x}\bigl(t^{\ast}\bigr) \bigr) \mathcal{Y}_{y;T}
\bigl(T-t^\ast\bigr) \bigr],\label{eq:forward-reverse-density}
\end{equation}
expressed with the variant of the reverse process introduced in \eqref{so1}
above.
\end{corollary}

%re3.6 #&#
\begin{remark}
\label{rem:cond-diff-pathdep} For fixed $x,y\in\mathbb{R}^{d}$ and
$s_{0}<t^{\ast}<T$ as before, let us define a process $Z$ by
\[
Z(t) \coloneqq Y_{y;T}(T-t), \qquad t^{\ast}\leq t\leq T.
\]
The idea is that we run along the reverse diffusion $Y$ backwards in time.
Then (\ref{cdy1}) reads
\begin{eqnarray*}
&&\mathbb{E} \bigl[  g \bigl( X_{s_{0},x}(s_{1}), \ldots,
X_{s_{0},x}\bigl(t^{\ast}\bigr), X_{s_{0},x}(t_{1}),
\ldots, X_{s_{0},x}(t_{L-1}) \bigr) \rrvert X_{s_{0},x}(T)
= y \bigr]
\\
&&\qquad=\frac{1}{p(s_{0},x,T,y)}\lim_{\epsilon\downarrow0}\mathbb{E}%
\bigl[g
\bigl(X_{s_{0},x}(s_{1}),\ldots,X_{s_{0},x}
\bigl(t^{\ast
}\bigr),Z(t_{1}),\ldots,Z(t_{L-1}%
)\bigr) \\
&&\hspace*{127pt}\qquad\quad{}\times K_{\epsilon} \bigl( Z\bigl(t^{\ast}\bigr)-X_{s_{0},x}
\bigl(t^{\ast}\bigr) \bigr) \mathcal{Y}_{y;T}
\bigl(T-t^\ast\bigr) \bigr].
\end{eqnarray*}
\end{remark}

%s3.2 #&#
\subsection{Representations for conditioning on a set}
\label{sec:repr-cond-set}

Now let us assume that we are interested in the conditional expectation
of a
functional $g ( X_{s_{0},x}(s_{1}),\break \ldots,  X_{s_{0},x}(t_{L-1})
) $
given $X_{T}\in A$ for some Borel set $A$. It is assumed for simplicity that
either $A$ is a subset of $\mathbb{R}^{d}$ with positive Lebesgue
measure and
with $\mathbb{P}(X_{s_{0},x}(T)\in A)>0$, or $A$ is an affine
hyperplane of
dimension $d^{\prime}$, $0\leq d^{\prime}\leq d$. As a further
simplification in
the latter case, although without further loss of generality, we assume that
$A$ is of the form
%
%e23 #&#
\begin{equation}
A= \bigl\{ x\in\mathbb{R}^{d}\dvtx x^{1}=c^{1},\ldots,x^{d-d^{\prime
}}=c^{d-d^{\prime
}}
\bigr\}. \label{hyp}
\end{equation}
For $0\leq d^{\prime}\leq d$ we consider the ``restricted'' Lebesgue measure
%
%e24 #&#
\begin{equation}
\lambda_{A}(dx)=\delta_{ \{ c^{1} \} }\bigl(dx^{1}\bigr)
\cdots\delta_{ \{ c^{d-d^{\prime}} \} }\bigl(dx^{d-d^{\prime
}}\bigr)\cdot
dx^{d-d^{\prime}+1}\cdots dx^{d}, \label{a}
\end{equation}
which coincides with the ordinary Lebesgue measure if $d^{\prime}=d$,
and with
a Dirac point measure if $d^{\prime}=0$. We next introduce a random variable
$\xi$ with support in $A$ independent from $X$ and $Y$, whose law has a density
$\varphi>0$ with respect to $\lambda_{A}$. Let further $(Y_{t^{\ast},\xi
},\mathcal{Y}_{\xi;T})$ denote the reverse process starting at the
random location $(\xi,1)$ at time $t^{\ast}$. Here, we replace
condition \eqref{assp} on the positivity of the transition density by
%
%e25 #&#
\begin{equation}
\label{eq:assp2} \int_{A}p(s_{0},x,T,z)
\lambda_{A}(dz)>0.
\end{equation}

%th3.7 #&#
\begin{theorem}
\label{thr:cond-dist-set} Let the kernel function $K$ be as in
Theorem \ref{thr:conditional-dist-general-grid}, and let there be given
a time
grid of the form \eqref{eq:full-grid}. The conditional expectation of
\[
g \bigl( X_{s_{0},x}(s_{1}),\ldots,X_{s_{0},x}
\bigl(t^{\ast}\bigr),X_{s_{0},x}(t_{1}), \ldots,
X_{s_{0},x}(t_{L-1}) \bigr)
\]
given $X_{s_{0},x}(T)\in A$ with $A$ being a Borel set, either with positive
probability or a hyperplane of the form (\ref{hyp}), and $g$ being a bounded
measurable test function, has the stochastic representation
\begin{eqnarray*}
&&\int_{A}p(s_{0},x,T,y)\lambda_{A}(dy)
\cdot\mathbb{E} \bigl[  g \bigl( X_{s_{0},x}(s_{1}),
\ldots, X_{s_{0},x}(t_{L-1}) \bigr) \rrvert X_{s_{0},x}(T)
\in A \bigr]
\\
&&\qquad=\lim_{\epsilon\downarrow0} \mathbb{E} \biggl[g \bigl( X_{s_{0},x}(s_{1}),
\ldots, X_{s_{0},x}\bigl(t^{\ast}\bigr), Y_{\xi;T}(
\widehat{t}_{L-1}), \ldots, Y_{\xi;T}(\widehat{t}_{1})
\bigr)
\\
&&\hspace*{96pt}\qquad\quad{}\times K_{\epsilon} \bigl( Y_{\xi;T}(\widehat{t}_L)-X_{s_{0},x}
\bigl(t^{\ast}\bigr) \bigr) \frac{\mathcal{Y}_{\xi;T}(\widehat{t}_L)}{\varphi(\xi)} \biggr].
\end{eqnarray*}
In particular, by setting $g\equiv1$ we obtain a stochastic
representation for
the factor
\[
\int_{A}p(s_{0},x,T,y)\lambda_{A}(dy)=
\lim_{\epsilon\downarrow0}%
\mathbb{E} \biggl[K_{\epsilon}
\bigl( Y_{\xi;T}(\widehat {t}_L)-X_{s_{0},x}
\bigl(t^{\ast
}\bigr) \bigr) \frac{\mathcal{Y}_{\xi;T}(\widehat{t}_L)}{\varphi(\xi)} \biggr].
\]
\end{theorem}
\begin{pf}
Let us abbreviate
\begin{eqnarray*}
H_{A} & \coloneqq&\mathbb{E} \bigl[  g \bigl(
X_{s_{0},x}(s_{1}%
),\ldots,X_{s_{0},x}(t_{L-1})
\bigr) \rrvert X_{s_{0},x}(T)\in A \bigr],
\\
H(y) & \coloneqq&\mathbb{E} \bigl[  g \bigl( X_{s_{0},x}(s_{1}),
\ldots,X_{s_{0},x}(t_{L-1}) \bigr) \rrvert X_{s_{0},x}(T)=y
\bigr],
\end{eqnarray*}
and consider the density of the conditional distribution of $X_{s_{0},x}(T)$
given $X_{s_{0},x}(T)\in A$ with respect to the measure $\lambda_{A}$,
that is,
\[
q(y)=\frac{p(s_{0},x,T,y)}{\int_{A}p(s_{0},x,T,z)\lambda_{A}(dz)}1_{A}(y).
\]
Recall \eqref{eq:assp2} and the construction (\ref{a}) of $\lambda
_{A}$. Then
we have
\begin{eqnarray*}
H_{A} &=& \int_{A}\mathbb{E} \bigl[  g
\bigl( X_{s_{0},x}%
(s_{1}),\ldots,X_{s_{0},x}(t_{L-1})
\bigr) \rrvert X_{s_{0}%
,x}(T)=y \bigr] q(y)\lambda_{A}(dy)
\\
&=&\int_{A}H(y)q(y)\lambda_{A}(dy)
\\
&=& \mathbb{E} \biggl[ \frac{H(\xi)q(\xi
)}{\varphi(\xi)} \biggr]
\\
&=& \frac{\mathbb{E} [
({p(s_{0},x,T,\xi)}/{\varphi(\xi)})\mathbb{E}  [  g (
X_{s_{0},x}(s_{1}),\ldots,X_{s_{0},x}(t_{L-1}) ) \rrvert
X_{s_{0},x}(T)=\xi ]  ]}{\int_{A}p(s_{0},x,T,z)\lambda_{A}(dz)}.
\end{eqnarray*}
Hence, denoting $p_A \coloneqq\int_{A}p(s_{0},x,T,z)\lambda_{A}(dz)$,
\begin{eqnarray*}
&& H_{A}\times p_A \\
&&\qquad= \mathbb{E} \biggl[
\frac{1}{\varphi(\xi)}\lim_{\epsilon\downarrow0} \mathbb{E}^{\xi} \bigl[g
\bigl( X_{s_{0},x}(s_{1}),\ldots,X_{s_{0},x}
\bigl(t^{\ast
}\bigr),Y_{\xi;T}(\widehat{t}_{L-1}),
\ldots,Y_{\xi;T}(\widehat{t}_{1}) \bigr)
\\
&&\hspace*{138pt}\qquad\quad{} \times K_{\epsilon} \bigl( Y_{\xi;T}(\widehat{t}_L)
- X_{s_{0},x}\bigl(t^{\ast}\bigr) \bigr) \mathcal{Y}_{\xi;T}(
\widehat{t}_L) \bigr] \biggr]
\\
&&\qquad=\lim_{\epsilon\downarrow0} \mathbb{E} \biggl[g \bigl( X_{s_{0},x}(s_{1}),
\ldots, X_{s_{0},x}\bigl(t^{\ast}\bigr), Y_{\xi;T}(
\widehat{t}_{L-1}), \ldots, Y_{\xi;T}(\widehat{t}_{1})
\bigr)
\\
&&\hspace*{95pt}\qquad\quad{} \times K_{\epsilon} \bigl( Y_{\xi;T}(\widehat{t}_L)
- X_{s_{0},x}\bigl(t^{\ast}\bigr) \bigr) \frac{\mathcal{Y}_{\xi;T}(\widehat{t}_L)}{\varphi(\xi)} \biggr].
\end{eqnarray*}
\upqed\end{pf}

%co3.8 #&#
\begin{corollary}
\label{cor:cond-dist-comp} The conditional expectation of $g (
X_{s_{0},x}(s_{1}),\ldots,\break  X_{s_{0},x}(t_{L-1}) ) $ given $X_{s_{0}%
,x}^{1}(T)=c^{1}\in\mathbb{R}$ has the stochastic representation
\begin{eqnarray*}
&&\lim_{\epsilon\downarrow0}\mathbb{E} \biggl[g \bigl( X_{s_{0},x}(s_{1}),
\ldots, X_{s_{0},x}\bigl(t^{\ast}\bigr),Y_{\xi;T}(
\widehat{t}_{L-1}),\ldots,Y_{\xi;T}(\widehat{t}_{1})
\bigr)
\\
&&\hspace*{84pt}\quad{}\times K_{\epsilon} \bigl( Y_{\xi;T}(\widehat{t}_L)-X_{s_{0},x}
\bigl(t^{\ast
}\bigr) \bigr) \frac{\mathcal{Y}_{\xi;T}(\widehat{t}_L)}{\varphi(\xi)} \biggr]\\
&&\qquad=
\int_{\mathbb{R}^{d-1}}p\bigl(s_{0},x,T,c^{1},y
\bigr)\,dy \\
&&\quad\qquad{}\times \mathbb{E} \bigl[  g \bigl( X_{s_{0},x}(s_{1}),
\ldots,X_{s_{0}%
,x}(t_{L-1}) \bigr) \rrvert X_{s_{0},x}^{1}(T)=c^{1}
\bigr]
\end{eqnarray*}
for any $\xi$ taking values in the hyperplane $A\coloneqq\{z\in\mathbb
{R}%
^{d} | z^{1}=c^{1}\}$ such that $\varphi> 0$ is the density of the law of
$\xi$ with respect to $\lambda_{A}$ defined accordingly. In particular, by
setting $g\equiv1$, we obtain a stochastic representation for the marginal
density
\begin{eqnarray*}
&&\lim_{\epsilon\downarrow0}\mathbb{E} \biggl[K_{\epsilon} \bigl(
Y_{\xi;T}(\widehat{t}_L) - X_{s_{0},x}
\bigl(t^{\ast}\bigr) \bigr)
\frac{\mathcal{Y}_{\xi;T}(\widehat{t}_L)}{\varphi(\xi)} \biggr] \\
&&\qquad=\int
_{\mathbb{R}^{d-1}}p\bigl(s_{0},x,T,c^{1},y^{2},
\ldots,y^{d}\bigr)\,dy^{2}\cdots dy^{d}.
\end{eqnarray*}
\end{corollary}

%re3.9 #&#
\begin{remark}
\label{rem:markov-chain}
Without doubt it is possible to construct analogous stochastic
representations for conditional Markov chains in the spirit
of \citet{MSS2}. The details, however, are considered beyond the scope
of the present
paper.
\end{remark}

%s4 #&#
\section{Forward-reverse estimators and their analysis}

\label{analysis} \label{sec:constr-analys-forw} The stochastic representations
for the conditional diffusion problem (\ref{cp}), derived in the previous
section, naturally lead to respective Monte Carlo estimators. In this section
we analyze the accuracy of these estimators, under the following
assumptions. First we need suitably regularity of the transition
densities of
both forward and reverse processes.

%co4.1 #&#
\begin{condition}
\label{ass:bound-density} We assume that the diffusion $X$ as well as the
reverse diffusion $Y$ (not including $\mathcal{Y}$) defined in \eqref{MC12}
have $C^{\infty}$ transition densities $p(t,x,s,y)$ and $q(t,x,s,y)$,
respectively. Moreover, for fixed $N\in\mathbb{N}$, there are constants
$m_{N}\in\mathbb{N}$, $\nu_{N}>0$, $\lambda_{N}>0$, $K_{N}>0$ and $C_{0}>0$
such that for any multi-indices
$\alpha,\beta\in\mathbb{N}_{0}^{d}$ with $|\alpha|+|\beta|\leq N$ we have
\[
\bigl\llvert \partial_{x}^{\alpha}\partial_{y}^{\beta
}p(t,x,s,y)
\bigr\rrvert \leq\frac{K_{N}}{(s-t)^{\nu_{N}}}\exp \biggl( -\lambda_{N}
\frac{|y-x|^{2}}{(1+C_{0}^{2})(s-t)} \biggr),
\]
uniformly for $(t,x,s,y)\in(0,T]\times\mathbb{R}^{d}\times\mathbb
{R}^{d}$ and
similarly for $q$.
\end{condition}

%re4.2 #&#
\begin{remark}
\label{rem:bound-density-restriction}
In fact, for the theorems formulated as below, we only need
Condition \ref{ass:bound-density} for $N = 2$. Higher order versions only
become necessary in the context of Remark \ref{rem:higher-order-kernel}.
\end{remark}

%re4.3 #&#
\begin{remark}\label{rem:ass_bound_density_autonomous}
By the results of \citet{KS}, Corollary~3.25, Condition \ref
{ass:bound-density}
is satisfied \emph{in the autonomous case} provided that (the vector fields
driving) the forward diffusion $X$ and $Y$ satisfy a uniform H\"ormander
condition, and $a$ and $\sigma$ are bounded, and $C^{\infty}$ bounded;
that is,
all the derivatives are bounded as well. We know of no similar study for
nonautonomous stochastic differential equations. Of course, the seminal
work by \citet{Ar} gives upper (and lower) Gaussian bounds for the transition
density of time-dependent, but uniformly elliptic stochastic differential
equations. Moreover, \citet{CM} prove the existence and smoothness of
transition densities for time-dependent SDEs under H\"ormander conditions.

In any case, an extension of the Kusuoka--Stroock result to the
time-inhomogeneous case seems entirely possible, in particular since
we do
not consider time-derivatives, for instance, by first considering the
case of
piecewise constant coefficients.
\end{remark}

%co4.4 #&#
\begin{condition}
\label{ass:kernel-order}
The kernel $K$ satisfies $\int_{\mathbb{R}^{d}}K(v)\,dv=1$
and\break  $\int_{\mathbb{R}^d} v K(v) \,dv = 0$. Moreover, it has lighter tails
than a Gaussian density in the sense that there are constants $C,\alpha
>0$ and
$\beta\geq0$ such that
\[
K(v)\leq C\exp \bigl( -\alpha|v|^{2+\beta} \bigr), \qquad v\in
\mathbb{R}^{d}.
\]
\end{condition}
In many applications, one would probably choose a compactly supported kernel,
which trivially satisfies the above tail-condition. Finally, we also introduce
some further assumptions put forth for convenience, which could be easily
relaxed.

%co4.5 #&#
\begin{condition}
\label{ass:convenience}
The functional $g\dvtx \mathbb{R}^{(K+L-1)\times d} \to\mathbb{R}$ together
with its gradient and its Hessian are bounded. Moreover, the
coefficient $c$
in \eqref{MC12} is bounded.
\end{condition}

%re4.6 #&#
\begin{remark}
\label{rem:convenience}
Condition \ref{ass:convenience} could be replaced by a requirement of
polynomial boundedness.
\end{remark}

%s4.1 #&#
\subsection{Forward-reverse estimators for conditioning on a fixed state}
\label{sec:forw-reverse-estim-state}

Let us consider
%
%e26 #&#
%e27 #&#
\begin{eqnarray}
\label{eq:def-heps}h_{\epsilon}&\coloneqq&\mathbb{E} \biggl[g \bigl(
X_{s_{0},x}(s_{1}), \ldots, X_{s_{0},x}
\bigl(t^{\ast}\bigr), Y_{y;T}(\widehat{t}_{L-1}),
\ldots, Y_{y;T}(\widehat{t}_{1}) \bigr)
\nonumber
\\[-8pt]
\\[-8pt]
\nonumber
&&\hspace*{38pt}\qquad{}\times\epsilon^{-d} K \biggl( \frac{Y_{y;T}(\widehat{t}_L)-X_{s_{0},x}(t^{\ast})}{\epsilon} \biggr)
\mathcal{Y}_{y;T}(\widehat{t}_L) \biggr],
\end{eqnarray}
which can---and will---be computed using Monte Carlo simulation. Here, we
recall the definition of $\widehat{t}_i = t_L - t_{L-i}$ given
in \eqref{eq:hat-grid}. By Theorem \ref{thr:conditional-dist-general-grid},
$h_\epsilon$ converges to
%
%e28 #&#
\begin{equation}
\label{eq:def-h} h \coloneqq p(s_0, x, T, y) \mathbb{E} \bigl[
 g \bigl( X_{s_0,x}(s_1), \ldots, X_{s_0,x}(t_{L-1})
\bigr) \rrvert X_{s_0,x}(T) = y \bigr].
\end{equation}

%th4.7 #&#
\begin{theorem}
\label{bias} Assuming
Conditions \ref{ass:bound-density}, \ref{ass:kernel-order}
and \ref{ass:convenience}, there are constants $C, \epsilon_{0} > 0$ such
that the bias of the approximation $h_{\epsilon}$ can be bounded by
\[
|h - h_{\epsilon}| \le C \epsilon^{2}, \qquad 0 < \epsilon<
\epsilon_{0}.
\]
\end{theorem}
\begin{pf}
Changing variables $y_{0} \to v \coloneqq\frac{y_{0} -
x_{K}}{\epsilon}$ in Theorem \ref{Gen1}, we arrive at
\begin{eqnarray*}
&& h_{\epsilon}= \int g(x_{1}, \ldots, x_{K},
y_{1}, \ldots, y_{L-1}) K(v)\\
&&\hspace*{35pt}{}\times \prod
_{i=1}^{K} p(s_{i-1}, x_{i-1},
s_{i}, x_{i})  p\bigl(t^{\ast}, x_{K}+\delta v, t_{1},
y_{1}\bigr) \\
&&\hspace*{35pt}{}\times\prod_{i=2}^{L}
p(t_{i-1}, y_{i-1}, t_{i}, y_{i})
\,dx_{1} \cdots dx_{K} \,dv \,dy_{1} \cdots
dy_{L-1}.
\end{eqnarray*}
In particular, we have that $h = \lim_{\epsilon\downarrow0} h_{\epsilon
}$. Consider
\[
r_{\epsilon}(x_{K}, y_{1}) \coloneqq\int k(v) p
\bigl(t^{\ast}, x_{K}+\epsilon v, t_{1},
y_{1}\bigr) \,dv - p\bigl(t^{\ast}, x_{K},
t_{1}, y_{1}\bigr).
\]
In the following, we use the notation $\partial_x^\beta
\coloneqq\partial_{x^1}^{\beta^1} \cdots\partial_{x^d}^{\beta^d}$,
for $x
\in\mathbb{R}^d$, $\beta\in\mathbb{N}^d$. By Taylor's formula,
Conditions \ref{ass:kernel-order} and~\ref{ass:bound-density},
we get
\begin{eqnarray*}
r_{\epsilon} & = &\int K(v) \bigl[ p\bigl(t^{\ast},
x_{K}+\epsilon v, t_{1}, y_{1}\bigr) - p
\bigl(t^{\ast}, x_{K}, t_{1}, y_{1}\bigr)
\bigr] \,dv
\\
& =& \int K(v) \bigl[ \epsilon\partial_{x} p\bigl(t^{\ast},
x_{K}, t_{1}, y_{1}\bigr) \cdot v \bigr] \,dv
\\
&&{} + \sum_{|\beta| = 2} \frac{2}{\beta!}
\epsilon^{2} \int\!\!\!\int_{0}^{1} (1-t)
\partial^{\beta}_{x} p\bigl(t^{\ast},x_{K} +
t \epsilon v, t_{1}, y_{1}\bigr) \cdot v^{\beta} \,dt
K(v) \,dv
\end{eqnarray*}
implying that
\begin{eqnarray*}
\bigl|r_{\epsilon}(x_{K}, y_{1})\bigr| & \le&\sum
_{|\beta| = 2} \frac{2}{\beta!} \epsilon^{2} \int
_{0}^{1} (1-t)\int\bigl\llvert
\partial^{\beta}_{x} p\bigl(t^{\ast
},x_{K} +
t \epsilon v, t_{1}, y_{1}\bigr) \bigr\rrvert
|v|^{2} K(v) \,dv \,dt
\\
& \le&\sum_{|\beta| = 2} \frac{2}{\beta!}
\epsilon^{2} C_{1} \int_{0}^{1}
(1-t) \int e^{-\gamma|y_{1} - x_{K} - t\epsilon v|^{2}} |v|^{2} K(v) \,dv \,dt
\\
& \le&\sum_{|\beta| = 2} \frac{2}{\beta!}
\epsilon^{2} C_{1} C_{\eta}\int
_{0}^{1} (1-t) \int e^{-\gamma|y_{1} - x_{K} - t\epsilon v|^{2}}
e^{-\eta|v|^{2}} \,dv \,dt,
\end{eqnarray*}
where\vspace*{-1pt} $C_{1} \coloneqq\frac{K_{2}(t_{1}-t^{\ast})}{(t_{1}-t^{\ast
})^{\nu_{2}}}$,
$\gamma\coloneqq\frac{\lambda_{2}}{(1+C_{0}^{2}) (t_{1} - t^{\ast})}$
as given in
Condition \ref{ass:bound-density}, $\eta> 0$ and $C_{\eta}$ is chosen such
that $|v|^{2} K(v) \le C_{\eta}e^{-\eta|v|^{2}}$, which is possible by
Condition \ref{ass:kernel-order}. Since
\[
|y_{1} - x_{K} - t\epsilon v|^{2} =
|y_{1} - x_{K}|^{2} - 2t\epsilon \langle
y_{1}-x_{K}, v \rangle+ t^{2}
\epsilon^{2} |v|^{2},
\]
we can further compute, using $\sigma^{2} \coloneqq\frac{1}{2 (
\eta+\gamma t^{2} \epsilon^{2}  ) } \le\frac{1}{2\eta}$,
\begin{eqnarray*}
&&\int e^{-\gamma|y_{1} - x_{K} - t\epsilon v|^{2}} e^{-\eta|v|^{2}} \,dv\\
&&\qquad = e^{-\gamma|y_{1}-x_{K}|^{2}} \int
e^{2t\gamma\epsilon\langle
y_{1}-x_{K},
v \rangle} e^{-{|v|^{2}}/{(2\sigma^{2})}} \,dv
\\
&&\qquad = \biggl( \frac{\eta+\gamma t^{2} \epsilon^{2}}{\pi} \biggr) ^{-d/2} \exp \biggl(
\epsilon^{2} \frac{t^{2} \gamma^{2}}{\eta} |y_{1}-x_{K}|^{2}
\biggr) e^{-\gamma|y_{1}-x_{K}|^{2}}
\\
&&\qquad \le \biggl( \frac{\pi}{\eta} \biggr) ^{d/2} e^{\epsilon^{2} ({\gamma
^{2}}/{\eta} )|y_{1} - x_{K}|^{2}}
e^{-\gamma|y_{1}-x_{K}|^{2}}.
\end{eqnarray*}
Defining $\widetilde{C}_{\eta}\coloneqq\sum_{|\beta|=2} \frac{1}{\beta
!} C_{1}
C_{\eta}(\pi/\eta)^{d/2}$, we get the bound
\begin{eqnarray*}
\bigl|r_{\epsilon}(x_{K}, y_{1})\bigr| & \le&2 \widetilde{C}
\epsilon^{2} e^{-\gamma|y_{1}-x_{K}|^{2}} \int_{0}^{1}
(1-t) e^{\epsilon^{2} ({\gamma^{2}}/{\eta}) |y_{1}-x_{K}|^{2}} \,dt
\\
& \le&\widetilde{C} \epsilon^{2} e^{-\gamma^{\prime}|y_{1}-x_{K}|^{2}},
\end{eqnarray*}
with $\gamma^{\prime} = \gamma- \frac{\gamma^{2}}{\eta} \epsilon^2$,
which is
positive for $0 < \epsilon< \epsilon_{0} \coloneqq
(\eta/\gamma)^{1/2}$. Consequently, for $0 < \epsilon< \epsilon_{0}$,
we can
interpret $s_{\epsilon}(x_{K},y_{1}) \coloneqq|r_{\epsilon}(x_{K},
y_{1})| /
(C_{2}\epsilon^{2})$ as a (Gaussian) transition density, which has
moments of
all orders, for a suitable normalization constant $C_{2}$,\ for which we can
derive explicit upper bounds. Thus we finally obtain
%
%e29 #&#
\begin{eqnarray} \label{eq:bias-bound-explicit}
|h_{\epsilon}- h| & \le&\int\bigl|g(x_{1}, \ldots, x_{K},
y_{1}, \ldots, y_{L-1})\bigr| \prod_{i=1}^{K}
p(s_{i-1}, x_{i-1}, s_{i}, x_{i})
\nonumber
\\
&&\hspace*{8pt}{} \times\bigl|r_{\epsilon}(x_{K}, y_{1})\bigr| \prod
_{i=2}^{L} p(t_{i-1}, y_{i-1},
t_{i}, y_{i}) \,dx_{1} \cdots dx_{K}
\,dy_{1} \cdots dy_{L-1}
\nonumber
\\
& \le &C_{2} \epsilon^{2} \int\bigl|g(x_{1}, \ldots,
x_{K}, y_{1}, \ldots, y_{L-1})\bigr| \prod
_{i=1}^{K} p(s_{i-1}, x_{i-1},
s_{i}, x_{i})
\\
&&\hspace*{32pt}{} \times s_{\epsilon}(x_{K}, y_{1}) \prod
_{i=2}^{L} p(t_{i-1}, y_{i-1},
t_{i}, y_{i}) \,dx_{1} \cdots dx_{K}
\,dy_{1} \cdots dy_{L-1}\nonumber
\\
& \eqqcolon& C \epsilon^{2} < \infty,
\nonumber
\end{eqnarray}
provided that $0 < \epsilon< \epsilon_{0}$, as the last expression can be
interpreted as
\[
C_{2} \mathbb{E} \bigl[ \bigl|g(Z_{s_{1}}, \ldots,
Z_{s_{K}}, Z_{t_{1}}, \ldots, Z_{t_{L-1}})\bigr| |
Z_{s_{0}} = x, Z_{T} = y \bigr]
\]
for a Markov process $Z$ with transition densities $p(s_{i-1}, x_{i-1}, s_{i},
x_{i})$, $1 \le i \le K$, $s_{\epsilon}(x_{K}, y_{1})$, $p(t_{i-1}, y_{i-1},
t_{i}, y_{i})$, $2 \le i \le L$, which admits finite moments of all orders
by construction.
\end{pf}

%re4.8 #&#
\begin{remark}
Note that the constant $C$ in the above statement can be explicitly bounded
in terms of the bound on $g$, the constants appearing in
Condition \ref{ass:bound-density} and $\eta$.
\end{remark}

In the spirit of \citet{MSS1} we now introduce a
Monte Carlo estimator $\widehat{h}_{\epsilon}$ for the quantity
$h_{\epsilon}$
introduced in \eqref{eq:def-heps}. Let us denote
%
%e30 #&#
%e31 #&#
\begin{eqnarray}
\label{eq:def-Z} Z_{nm}^\epsilon&\coloneqq&\frac{1}{\epsilon^{d}} g
\bigl( X_{s_{0},x}^{n}(s_{1}), \ldots,
X_{s_{0},x}^{n}(s_{K}), Y_{y;T}^{m}(
\widehat{t}_{L-1}), \ldots, Y_{y;T}^{m}(
\widehat{t}_{1}) \bigr)
\nonumber
\\[-8pt]
\\[-8pt]
\nonumber
&&{}\times K \biggl( \frac{Y_{y;T}^{m}(\widehat{t}_L) - X_{s_{0},x}^{n}(t^{\ast})}{\epsilon
} \biggr) \mathcal{Y}_{y;T}^{m}(
\widehat{t}_L).
\end{eqnarray}
Note that $\mathbb{E} [ Z^\epsilon_{nm} ] = h_{\epsilon}$.
The Monte Carlo estimator is now defined by
%
%e32 #&#
\begin{equation}
\label{eq:hat-heps} \widehat{h}_{\epsilon, M, N} \coloneqq\frac{1}{NM} \sum
_{n=1}^N \sum
_{m=1}^M Z_{nm}^\epsilon,
\end{equation}
where the superscripts $n$ and $m$ denote different, independent realizations
of the corresponding processes.\vadjust{\goodbreak} We are left to analyze the variance of the
estimator $\widehat{h}_{\epsilon,M,N}$. To this end, we consider the
expectation $\mathbb{E} [ Z^\epsilon_{nm} Z^\epsilon_{n^{\prime
}m^{\prime}} ] $ for various combinations of $n,m,n^{\prime}$ and
$m^{\prime}$.

%re4.9 #&#
\begin{remark}
For the remainder of the section, we omit the sub-scripts in $X$, $Y$ and
$\mathcal{Y}$ as we keep the initial times and values fixed.
\end{remark}

%le4.10 #&#
\begin{lemma}
\label{lem:Z-m-mprime} For $m \neq m^{\prime}$ we obtain
\begin{eqnarray*}
&&\mathbb{E} \bigl[Z^\epsilon_{nm}Z^\epsilon_{nm^{\prime}}
\bigr] \rrvert _{\epsilon=0} \\
&&\qquad= \int g(x_{1}, \ldots,
x_{K}, y_{1}, \ldots, y_{L-1}) g
\bigl(x_{1}, \ldots, x_{K}, y_{1}^{\prime},
\ldots, y_{L-1}^{\prime}\bigr)
\\
&&\hspace*{9pt}\qquad\quad{}\times p\bigl(t^{\ast}, x_{K}, t_{1},
y_{1}\bigr) p\bigl(t^{\ast}, x_{K},
t_{1}, y_{1}^{\prime}\bigr)
\\
&&\hspace*{9pt}\qquad\quad{}\times\prod_{i=1}^{K}
p(s_{i-1}, x_{i-1}, s_{i}, x_{i})
\,dx_{i} \prod_{i=2}^{L}
p(t_{i-1}, y_{i-1}, t_{i}, y_{i})
\,dy_{i-1} \\
&&\hspace*{9pt}\quad\qquad{}\times\prod_{i=2}^{L} p
\bigl(t_{i-1}, y_{i-1}^{\prime}, t_{i},
y_{i}^{\prime}\bigr) \,dy_{i-1}^{\prime}.
\end{eqnarray*}
Moreover, we can bound
\[
\bigl\llvert \mathbb{E}\bigl[Z^\epsilon_{nm}
Z^\epsilon_{nm^{\prime}}\bigr] - \mathbb{E} \bigl[Z^\epsilon_{nm}
Z^\epsilon_{nm^{\prime}}\bigr] \rrvert _{\epsilon=0} \bigr\rrvert \le
C \epsilon^{2}.
\]
\end{lemma}

\begin{pf}
In what follows, $C$ is a positive constant, which may change from line to
line. We have
\begin{eqnarray*}
&&E \bigl[ Z^\epsilon_{nm} Z^\epsilon_{nm^{\prime}}
\bigr] \\
&&\qquad= \epsilon^{-2d} E \biggl[ g \bigl( X_{s_{1}}^{n},
\ldots, X_{s_{K}}^{n}, Y_{\widehat{t}_{L-1}}^{m},
\ldots, Y_{\widehat{t}_{1}}^{m} \bigr) g \bigl( X_{s_{1}}^{n},
\ldots, X_{s_{K}}^{n}, Y_{\widehat{t}_{L-1}}^{m^{\prime}},
\ldots, Y_{\widehat{t}_{1}}^{m^{\prime}} \bigr)
\\
&&\hspace*{129pt}\qquad\quad{} \times K \biggl( \frac{Y_{\widehat{t}_L}^{m}-X_{t^{\ast}}^{n}}{\epsilon
} \biggr) K \biggl( \frac{Y_{\widehat{t}_L}^{m^{\prime}}-X_{t^{\ast}}^{n}}{\epsilon
}
\biggr) \mathcal{Y}_{\widehat{t}_L}^{m} \mathcal{Y}_{\widehat{t}_L}^{m^{\prime
}}
\biggr]
\\
&&\qquad = \epsilon^{-2d} \int g(x_{1}, \ldots, x_{K},
y_{1}, \ldots, y_{L-1}) g\bigl(x_{1}, \ldots,
x_{K}, y_{1}^{\prime}, \ldots, y_{L-1}^{\prime}
\bigr)
\\
& &\hspace*{32pt}\qquad\quad{}\times K \biggl( \frac{y_{0}-x_{K}}{\epsilon} \biggr) K \biggl( \frac{y_{0}^{\prime}-x_{K}}{\epsilon}
\biggr) \prod_{i=1}^{K}
p(s_{i-1}, x_{i-1}, s_{i}, x_{i})
\,dx_{i}
\\
& &\hspace*{32pt}\qquad\quad{}\times\prod_{i=1}^{L}
p(t_{i-1}, y_{i-1}, t_{i}, y_{i})
\,dy_{i-1} \prod_{i=1}^{L} p
\bigl(t_{i-1}, y_{i-1}^{\prime}, t_{i},
y_{i}^{\prime}\bigr) \,dy_{i-1}^{\prime}
\\
&&\qquad = \int g(x_{1}, \ldots, x_{K}, y_{1}, \ldots,
y_{L-1}) g\bigl(x_{1}, \ldots, x_{K},
y_{1}^{\prime}, \ldots, y_{L-1}^{\prime}\bigr)
\\
&&\hspace*{8pt}\qquad\quad{} \times K(v) K\bigl(v^{\prime}\bigr) p\bigl(t^\ast,
x_{K}+\epsilon v, t_{1}, y_{1}\bigr) \,dv p
\bigl(t^{\ast}, x_{K} + \epsilon v^{\prime},
t_{1}, y_{1}^{\prime}\bigr) \,dv^{\prime
}
\\
&&\hspace*{8pt}\qquad\quad{} \times\prod_{i=1}^{K}
p(s_{i-1}, x_{i-1}, s_{i}, x_{i})
\,dx_{i}\\
&&\hspace*{8pt}\qquad\quad{} \times \prod_{i=2}^{L}
p(t_{i-1}, y_{i-1}, t_{i}, y_{i})
\,dy_{i-1}
\\
&&\hspace*{8pt}\qquad\quad{} \times\prod_{i=2}^{L} p
\bigl(t_{i-1}, y_{i-1}^{\prime}, t_{i},
y_{i}^{\prime}\bigr) \,dy_{i-1}^{\prime},
\end{eqnarray*}
where we changed variables $v \coloneqq(y_{0}-x_{K})/\epsilon$ and
$v^{\prime}\coloneqq(y_{0}^{\prime}-x_{K})/\epsilon$. Thus, for
$\epsilon=0$,
we arrive at the above expression, which is treated as a problem-dependent
constant.

Using Condition \ref{ass:kernel-order} [and the short-hand notation $p(x,y)
\coloneqq p(t^\ast,x,t_1, y)$], we now consider
\begin{eqnarray*}
% r_{\epsilon}^{(1,2)}(x_{K},y_{1},y_{1}^{\prime}) &\coloneqq\int K(v)
% K(v^{\prime}) \bigl[ p(t^{\ast}, x_{K}+\epsilon v, t_{1}, y_{1}) p(t^{
% x_{K}+\epsilon v^{\prime}, t_{1}, y_{1}^{\prime}) -\\
% & - p(t^{\ast}, x_{K}, t_{1}, y_{1}) p(t^{\ast}, x_{K}, t_{1},
% y_{1}^{\prime}) \Bigr] dv dv^{\prime}\\
% & = \epsilon^{2} \int_{0}^{1} (1-t) \Biggl[ \sum_{i=1}^{d} \int K(v)
% K(v^{\prime}) \partial_{x}^{2e_{i}}p(t^{\ast}, x_{K}+t\epsilon v,
%t_{1},
% y_{1}) p(t^{\ast}, x_{K}+t\epsilon v^{\prime}, t_{1}, y_{1}^{\prime})
% v_{i}^{2} dv dv^{\prime}+\\
% & + \sum_{i=1}^{d} \int K(v) K(v^{\prime\ast}, x_{K}+t\epsilon v,
%t_{1},
% y_{1}) \partial_{x}^{2e_{i}} p(t^{\ast}, x_{K}+t\epsilon v^{\prime},
%t_{1},
% y_{1}^{\prime}) (v^{\prime}_{i})^{2} dv dv^{\prime}+\\
% & + 2 \sum_{i,j=1}^{d} \int K(v) K(v^{\prime}) \partial_{x}^{e_{i}}
% p(t^{\ast}, x_{K}+t\epsilon v, t_{1}, y_{1}) \partial_{x}^{e_{j}}
%p(t^{\ast},
% x_{K}+t\epsilon v^{\prime}, t_{1}, y_{1}^{\prime}) v_{i} v^{
% dv^{\prime}\Biggr] dt,
&&\hspace*{-4pt}r_{\epsilon}^{(1,2)}
\bigl(x_{K},y_{1},y_{1}^{\prime}\bigr) \\
&&\hspace*{-6pt}\qquad
\coloneqq \int K(v) K\bigl(v^{\prime}\bigr) \\
&&\hspace*{-6pt}\hspace*{44pt}{}\times\bigl[ p(x_{K}+
\epsilon v, y_{1}) p\bigl(x_{K}+\epsilon v^{\prime},
y_{1}^{\prime}\bigr)
 - p(x_{K}, y_{1}) p\bigl(x_{K},
y_{1}^{\prime}\bigr) \bigr] \,dv \,dv^{\prime}
\\
&&\hspace*{-6pt}\qquad = \epsilon^{2} \int K(v) K\bigl(v^\prime\bigr)\\
&&\hspace*{-6pt}\hspace*{21pt}\qquad\quad{}\times \int
_{0}^{1} (1-t) \Biggl[ \sum
_{i=1}^{d} \partial_{x}^{2e_{i}}p(x_{K}+t
\epsilon v, y_{1}) p\bigl(x_{K}+t\epsilon v^{\prime},
y_{1}^{\prime}\bigr) v_{i}^{2}
\\
&&\hspace*{-6pt}\hspace*{84pt}\qquad\quad{} + \sum_{i=1}^{d} p(x_{K}+t
\epsilon v, y_{1}) \\
&&\hspace*{-6pt}\hspace*{122pt}\qquad{}\times\partial_{x}^{2e_{i}} p
\bigl(x_{K}+t\epsilon v^{\prime}, y_{1}^{\prime}
\bigr) \bigl(v^{\prime}_{i}\bigr)^{2}
\\
&&\hspace*{-6pt}\hspace*{-6pt}\hspace*{84pt}\qquad\quad{} + 2 \sum_{i,j=1}^{d}
\partial_{x}^{e_{i}} p(x_{K}+t\epsilon v,
y_{1}) \\
&&\hspace*{-6pt}\hspace*{156pt}{}\times\partial_{x}^{e_{j}} p\bigl(x_{K}+t
\epsilon v^{\prime}, y_{1}^{\prime}\bigr) v_{i}
v^{\prime}_{j} \,dv \,dv^{\prime} \Biggr] \,dt \,dv
\,dv^{\prime},
\end{eqnarray*}
where, for instance, $\partial_x^{e_i} \equiv\partial_{x^i}$ and
$\partial_x^{2e_i} \equiv\partial_{x^i} \partial_{x^i}$. By similar
techniques as in the proof of Theorem \ref{bias}, relying once more on the
uniform bounds of Condition \ref{ass:bound-density}, we arrive at an upper
bound
\[
\bigl\llvert r_{\epsilon}^{(1,2)}\bigl(x_{K},y_{1},y_{1}^{\prime}
\bigr) \bigr\rrvert \le C s^{(1,2)}_{\epsilon}(x_{K},
y_{1}) s^{(1,2)}_{\epsilon}\bigl(x_{K},
y_{1}^{\prime
}\bigr)
\]
for a transition density $s^{(1,2)}_{\epsilon}(x_{K}, y_{1})$ with Gaussian
bounds. Consequently, we obtain
\begin{eqnarray*}
&&\bigl\llvert \mathbb{E}\bigl[Z^\epsilon_{nm}
Z^\epsilon_{nm^{\prime}}\bigr] - \mathbb{E} \bigl[Z^\epsilon_{nm}
Z^\epsilon_{nm^{\prime}}\bigr] \rrvert _{\epsilon=0} \bigr\rrvert \\
&&\qquad\le
C \epsilon^{2} \int\bigl|g(x_{1}, \ldots, x_{K},
y_{1}, \ldots, y_{L-1})\bigr|
\\
&&\hspace*{30pt}\qquad\quad{}\times\bigl|g\bigl(x_{1}, \ldots, x_{K}, y_{1}^{\prime},
\ldots, y_{L-1}^{\prime
}\bigr)\bigr|\\
 &&\hspace*{30pt}\qquad\quad{}\times\prod
_{i=1}^{K} p(s_{i-1}, x_{i-1},
s_{i}, x_{i}) \,dx_{i} s_{\epsilon}^{(1,2)}(x_{K},
y_{1}) \,dy_{1}
\\
&&\hspace*{30pt}\qquad\quad{}\times\prod_{i=3}^{L}
p(t_{i-1}, y_{i-1}, t_{i}, y_{i})
\,dy_{i-1} \times s_{\epsilon}^{(1,2)}\bigl(x_{K},
y_{1}^{\prime}\bigr) \,dy_{1}^{\prime}\\
&&\hspace*{30pt}\qquad\quad{}\times\prod
_{i=3}^{L} p\bigl(t_{i-1},
y_{i-1}^{\prime}, t_{i}, y_{i}^{\prime}
\bigr) \,dy_{i-1}^{\prime},
\end{eqnarray*}
which can be bounded by $C \epsilon^{2}$ by boundedness of $g$. In
fact, we
can find densities $\widetilde{p}$ and $\widetilde{q}$ with Gaussian tails
such that
%
%e33 #&#
\begin{eqnarray}
\label{eq:Zm-mprime}&&\bigl\llvert \mathbb{E}\bigl[Z^\epsilon_{nm}
Z^\epsilon _{nm^{\prime}}\bigr] - \mathbb{E}
\bigl[Z^\epsilon_{nm} Z^\epsilon_{nm^{\prime}}\bigr]
\rrvert _{\epsilon=0} \bigr\rrvert
\nonumber
\\[-8pt]
\\[-8pt]
\nonumber
&&\qquad \le C \epsilon^{2} \int
\widetilde{p}\bigl(s_{0},x,t^{\ast},x_{K}\bigr)
\widetilde{q}\bigl(t^{\ast},x_{K},T,y\bigr)^{2}
\,dx_{K}.%\qedhere
\end{eqnarray}
\upqed\end{pf}

When we consider $\mathbb{E} [ Z^\epsilon_{nm}
Z^\epsilon_{n^{\prime}m} ] $, we have to take care of terms
$\mathcal{Y}_{\widehat{t}_L}^{2}$ appearing in the expectation. To this
end, let us
introduce
\begin{eqnarray*}
\mu_{2}(y_{0}, \ldots, y_{L-1}) & \coloneqq
\mathbb{E} \bigl[  \mathcal{Y}_{\widehat{t}_L}^{2} \rrvert
Y_{\widehat{t}_L} = y_{0}, \ldots, Y_{\widehat{t}_{1}} =
y_{L-1} \bigr].
\end{eqnarray*}
In what follows, we replace $\mathcal{Y}_{\widehat{t}_L}^{2}$ by its
conditional expectation $\mu_{2} ( Y_{\widehat{t}_{L}}, \ldots,
Y_{\widehat{t}_{1}} )$ and re-write the expectation as an integral
w.r.t. the transition density of the reverse diffusion $Y$; by independence
of $X$ and $(Y,\mathcal{Y})$, we do not need to condition on $X$ as
well. Note
that by Condition \ref{ass:convenience}, $\mu_{2}$ is a bounded
function, and
the transition densities $q$ of the reverse process $Y$ satisfy the bounds
provided by Condition \ref{ass:bound-density} as well.

%le4.11 #&#
\begin{lemma}
\label{lem:Zn-nprime} For $n \neq n^{\prime}$ we have
\begin{eqnarray*}
 &&\mathbb{E} \bigl[ Z^\epsilon_{nm} Z^\epsilon_{n^{\prime}m}
\bigr] \rrvert _{\epsilon=0} \\
&&\qquad= \int g ( x_{1}, \ldots,
x_{K-1}, y_{0}, \ldots, y_{L-1} ) g \bigl(
x_{1}^{\prime}, \ldots, x_{K-1}^{\prime},
y_{0}, \ldots, y_{L-1} \bigr)
\\
&&\hspace*{9pt}\qquad\quad{}\times\mu_{2}(y_{0}, \ldots, y_{L-1}) \prod
_{i=1}^{K-1} p(s_{i-1},
x_{i-1}, s_{i}, x_{i}) \,dx_{i} \\
&&\hspace*{9pt}\qquad\quad{}\times{}\prod
_{i=1}^{K-1} p\bigl(s_{i-1},
x_{i-1}^{\prime}, s_{i}, x_{i}^{\prime}
\bigr) \,dx_{i}^{\prime}\times
\\
&&\hspace*{9pt}\qquad\quad{}\times p(s_{K-1}, x_{K-1}, s_{K},
y_{0}) p\bigl(s_{K-1}, x_{K-1}^{\prime},
s_{K}, y_{0}\bigr)\\
&&\hspace*{9pt}\qquad\quad{}\times \prod_{i=1}^{L}
q(\widehat{t}_{i-1}, y_{i}, \widehat{t}_{i},
y_{i-1}) \,dy_{i-1}.
\end{eqnarray*}
Moreover, there is a constant $C$ such that
\[
\bigl\llvert \mathbb{E}\bigl[Z^\epsilon_{nm}
Z^\epsilon_{n^{\prime}m}\bigr] - \mathbb{E} \bigl[Z^\epsilon_{nm}
Z^\epsilon_{n^{\prime}m}\bigr]\rrvert _{\epsilon=0}  \bigr\rrvert \le
\epsilon^{2} C.
\]
\end{lemma}

\begin{pf}
We first note that
\begin{eqnarray*}
&&\mathbb{E} \bigl[ Z_{nm}^\epsilon Z^\epsilon_{n'm}
\bigr] \\
&&\quad= \epsilon^{-2d} \mathbb{E} \biggl[ g \bigl(X^n_{s_1},
\ldots, X^n_{s_K}, Y^m_{\widehat{t}_{L-1}}, \ldots,
Y^m_{\widehat{t}_1} \bigr) g \bigl(X^{n'}_{s_1},
\ldots, X^{n'}_{s_K}, Y^m_{\widehat{t}_{L-1}}, \ldots,
Y^m_{\widehat{t}_1} \bigr)
\\
&&\hspace*{133pt}\qquad{} \times K \biggl( \f{Y^m_{\widehat{t}_L} - X^n_{t^\ast}}
{\epsilon} \biggr) K \biggl( \f{Y^m_{\widehat{t}_L} -
X^{n'}_{t^\ast}} {\epsilon} \biggr) \bigl( \mathcal{Y}^m_{\widehat{t}_L}
\bigr)^2 \biggr]
\\
&&\quad= \epsilon^{-2d} \mathbb{E} \biggl[ g \bigl(X^n_{s_1},
\ldots, X^n_{s_K}, Y^m_{\widehat{t}_{L-1}}, \ldots,
Y^m_{\widehat{t}_1} \bigr) g \bigl(X^{n'}_{s_1},
\ldots, X^{n'}_{s_K}, Y^m_{\widehat{t}_{L-1}}, \ldots,
Y^m_{\widehat{t}_1} \bigr)
\\
&&\hspace*{87pt}\qquad{} \times K \bigl( \f{Y^m_{\widehat{t}_L} - X^n_{t^\ast}}
{\epsilon} \bigr) K \biggl( \f{Y^m_{\widehat{t}_L} -
X^{n'}_{t^\ast}} {\epsilon} \biggr) \mu_2 \biggl(
Y^m_{\widehat{t}_L}, \ldots, Y^m_{\widehat{t}_1} \biggr)
\biggr].
\end{eqnarray*}

By a similar approach as in Lemma \ref{lem:Z-m-mprime}, but changing variables
$x_{K} \to v \coloneqq(y_{0}-x_{K})/\epsilon$ and $x_{K}^{\prime}\to
v^{\prime}\coloneqq(y_{0}-x_{K}^{\prime})/\epsilon$, we arrive at
\begin{eqnarray*}
&&\mathbb{E} \bigl[ Z^\epsilon_{nm} Z^\epsilon_{n^{\prime}m}
\bigr]\\
&&\qquad = \int g ( x_{1}, \ldots, x_{K-1}, y_{0}-
\epsilon v, y_{1}, \ldots, y_{L-1} )
\\
&&\hspace*{8pt}\qquad\quad{}\times g \bigl( x_{1}^{\prime}, \ldots, x_{K-1}^{\prime},
y_{0}-\epsilon v^{\prime}, y_{1}, \ldots,
y_{L-1} \bigr) \\
&&\hspace*{8pt}\qquad\quad{}\times  K ( v ) K \bigl( v^{\prime} \bigr)
\mu_{2}(y_{0}, \ldots, y_{L-1})
\\
&&\hspace*{8pt}\qquad\quad{}\times\prod_{i=1}^{K-1}
p(s_{i-1}, x_{i-1}, s_{i}, x_{i})
\,dx_{i} \prod_{i=1}^{K-1} p
\bigl(s_{i-1}, x_{i-1}^{\prime}, s_{i},
x_{i}^{\prime}\bigr) \,dx_{i}^{\prime}
\\
&&\hspace*{8pt}\qquad\quad{}\times p(s_{K-1}, x_{K-1}, s_{K},
y_{0}-\epsilon v) \,dv \\
&&\hspace*{8pt}\qquad\quad{}\times p\bigl(s_{K-1}, x_{K-1}^{\prime},
s_{K}, y_{0}-\epsilon v^{\prime}\bigr)
\,dv^{\prime}\\
&&\hspace*{8pt}\qquad\quad{}\times\prod_{i=1}^{L} q(
\widehat{t}_{i-1}, y_{i}, \widehat{t}_{i},
y_{i-1}) \,dy_{i-1}.
\end{eqnarray*}
For $\epsilon=0$, Condition \ref{ass:kernel-order} implies
\begin{eqnarray*}
&& \mathbb{E} \bigl[ Z^\epsilon_{nm} Z^\epsilon_{n^{\prime}m}
\bigr] \rrvert _{\epsilon=0} \\
&&\qquad= \int g ( x_{1}, \ldots,
x_{K-1}, y_{0}, \ldots, y_{L-1} ) g \bigl(
x_{1}^{\prime}, \ldots, x_{K-1}^{\prime},
y_{0}, \ldots, y_{L-1} \bigr)
\\
&&\hspace*{8pt}\qquad\quad{}\times\mu_{2}(y_{0}, \ldots, y_{L-1}) \prod
_{i=1}^{K-1} p(s_{i-1},
x_{i-1}, s_{i}, x_{i}) \,dx_{i} \\
&&\hspace*{8pt}\qquad\quad{}\times\prod
_{i=1}^{K-1} p\bigl(s_{i-1},
x_{i-1}^{\prime}, s_{i}, x_{i}^{\prime}
\bigr) \,dx_{i}^{\prime}
\\
&&\hspace*{8pt}\qquad\quad{}\times p(s_{K-1}, x_{K-1}, s_{K},
y_{0}) p\bigl(s_{K-1}, x_{K-1}^{\prime},
s_{K}, y_{0}\bigr) \\
&&\hspace*{8pt}\qquad\quad{}\times\prod_{i=1}^{L}
q(\widehat{t}_{i-1}, y_{i}, \widehat{t}_{i},
y_{i-1}) \,dy_{i-1},
\end{eqnarray*}
which gives the formula from the statement of the lemma.

For the bound on the difference, note once again that
\begin{eqnarray*}
r^{(2,1)}_{\epsilon}&\coloneqq&\int \bigl[ g(x_{1}, \ldots,
x_{K-1}, y_{0}-\epsilon v, y_{1}, \ldots,
y_{L-1}) \\
&&\hspace*{10pt}{}\times g\bigl(x_{1}^{\prime}, \ldots,
x_{K-1}^{\prime}, y_{0}-\epsilon v^{\prime},
y_{1}, \ldots, y_{L-1}\bigr)
\\
&&\hspace*{10pt}{}\times p(s_{K-1}, x_{K-1}, s_{K},
y_{0}-\epsilon v) \\
&&\hspace*{10pt}{}\times p\bigl(s_{K-1}, x_{K-1}^{\prime},
s_{K}, y_{0}-\epsilon v^{\prime}\bigr)
\\
&&\hspace*{10pt}{}- g(x_{1}, \ldots, x_{K-1}, y_{0}, \ldots,
y_{L-1})\\
&&\hspace*{10pt}{} \times g\bigl(x_{1}^{\prime}, \ldots,
x_{K-1}^{\prime}, y_{0}, \ldots, y_{L-1}
\bigr)
\\
&&\hspace*{10pt}{}\times p(s_{K-1}, x_{K-1}, s_{K}, y_{0}) p
\bigl(s_{K-1}, x_{K-1}^{\prime}, s_{K},
y_{0}\bigr) \bigr] K(v) K\bigl(v^{\prime}\bigr) \,dv
\,dv^{\prime}
\end{eqnarray*}
can be bounded in the sense that $|r^{(2,1)}_{\epsilon}| \le C
s^{(2,1)}_{\epsilon}(x_{K-1}, y_{0}) s^{(2,1)}_{\epsilon
}(x_{K-1}^{\prime},
y_{0})$ for transition densities $s^{(2,1)}_{\epsilon}$ with Gaussian tails,
so that
\begin{eqnarray*}
&&\bigl\llvert \mathbb{E}\bigl[Z^\epsilon_{nm}
Z^\epsilon_{n^{\prime}m}\bigr] - \mathbb{E} \bigl[Z^\epsilon_{nm}
Z^\epsilon_{n^{\prime}m}\bigr] \rrvert _{\epsilon=0} \bigr\rrvert\\
&&\qquad \le
C \epsilon^{2} \int\mu_{2}(y_{0}, \ldots,
y_{L-1}) \prod_{i=1}^{K-1}
p(s_{i-1}, x_{i-1}, s_{i}, x_{i})
\,dx_{i}
\\
&&\qquad\quad{}\times\prod_{i=1}^{K-1} p
\bigl(s_{i-1}, x_{i-1}^{\prime}, s_{i},
x_{i}^{\prime}\bigr) \,dx_{i}^{\prime}
s_{\epsilon}^{(2,1)}(x_{K-1}, y_{0})
s_{\epsilon}^{(2,1)}\bigl(x_{K-1}^{\prime},y_{0}
\bigr) \\
&&\qquad\quad{}\times \prod_{i=1}^{L} q(\widehat
{t}_{i-1}, y_{i}, \widehat{t}_{i},
y_{i-1}) \,dy_{i-1}.
\end{eqnarray*}
If $q$ was symmetric, that is, $q(\widehat{t}_{i-1}, y_{i}, \widehat{t}_{i},
y_{i-1}) = q(\widehat{t}_{i-1}, y_{i-1}, \widehat{t}_{i}, y_{i})$, then this
expression would already have the desired form. While symmetry of $q$
would be
a very strong assumption, note that Condition \ref{ass:bound-density}
allows us to bound
\[
q(\widehat{t}_{i-1}, y_{i}, \widehat{t}_{i},
y_{i-1}) \le C^{i} \exp \bigl( -\gamma^{i}
|y_{i} - y_{i-1}|^{2} \bigr) \eqqcolon
\widetilde{C}^is_{i-1}(y_{i-1}, y_{i})
\]
by a Gaussian transition density $s_{i-1}$ which is naturally symmetric.
Absorbing $\Vert\mu_{2} \Vert_{\infty}$ and $\prod_{i=1}^{L}
\widetilde{C}^{i}$ into the constant $C$ and denoting (by a mild abuse of
notation)
\begin{eqnarray*}
\widetilde{p}\bigl(s_{0},x,t^{\ast},y_{1}\bigr) &
\coloneqq&\int \prod_{i=1}^{K-1}
p(s_{i-1}, x_{i-1}, s_{i}, x_{i})
\,dx_{i} s_{\epsilon}^{(2,1)}(x_{K-1},y_{0}),
\\
\widetilde{q}\bigl(t^{\ast}, y_{0}, T, y\bigr) &\coloneqq&
\int\prod_{i=1}^{L} s_{i-1}(y_{i-1},
y_{i}) \,dy_{1} \cdots dy_{L-1},
\end{eqnarray*}
the Chapman--Kolmogorov equation implies that
%
%e34 #&#
\begin{eqnarray}
\label{eq:Zn-nprime} &&\bigl\llvert \mathbb{E}\bigl[Z^\epsilon_{nm}
Z^\epsilon_{n^{\prime}m}\bigr] - \mathbb{E} \bigl[Z^\epsilon_{nm}
Z^\epsilon_{n^{\prime}m}\bigr] \rrvert _{\epsilon=0} \bigr\rrvert \nonumber\\
&&\qquad
\le C \epsilon^{2} \int\widetilde{p}\bigl(s_{0},x,t^{\ast},y_{0}
\bigr)^{2} \widetilde{q}\bigl(t^{\ast}, y_{0}, T, y
\bigr) \,dy_{0}
\\
&&\qquad\le C \epsilon^{2} \int\widetilde{p}\bigl(s_{0},x,t^{\ast},y_{0}
\bigr) \widetilde {q}\bigl(t^{\ast}, y_{0}, T, y\bigr)
\,dy_{0}.\nonumber %%\qedhere
\end{eqnarray}
\upqed\end{pf}

%le4.12 #&#
\begin{lemma}
\label{lem:Znm} We have
\begin{eqnarray*}
&&\epsilon^{d} E \bigl[ \bigl(Z^\epsilon_{nm}
\bigr)^{2} \bigr]\\
&&\qquad = \int K(v)^{2} \,dv \int
g(x_{1}, \ldots, x_{K-1}, y_{0},
y_{1}, \ldots, y_{L-1})
\\
&&\hspace*{63pt}\qquad\quad{}\times\mu_{2}(y_{0}, y_{1}, \ldots,
y_{L-1}) \\
&&\hspace*{63pt}\qquad\quad{}\times\prod_{i=1}^{K-1}
p(s_{i-1}, x_{i-1}, s_{i}, x_{i})
p(s_{K-1}, x_{K-1}, s_{K}, y_{0})
\\
&&\hspace*{63pt}\qquad\quad{}\times\prod_{i=1}^{L} q(
\widehat{t}_{i-1}, y_{i}, \widehat{t}_{i},
y_{i-1}) \,dx_{1} \cdots dx_{K-1} \,dy_{0}
\,dy_{1} \cdots dy_{L-1}.
\end{eqnarray*}
Moreover, there is a constant $C > 0$ such that
\[
\Bigl\llvert \epsilon^{d} \mathbb{E}\bigl[\bigl(Z^\epsilon_{nm}
\bigr)^{2}\bigr] - \lim_{\epsilon
\to0} \epsilon^{d}
\mathbb{E} \bigl[ \bigl(Z^\epsilon_{nm}\bigr)^{2}
\bigr] \Bigr\rrvert \le C \epsilon^{2}.
\]
\end{lemma}

\begin{pf}
Substituting $x_{K} \to v \coloneqq(y_{0}-x_{K})/\epsilon$, we obtain
\begin{eqnarray*}
&&\epsilon^{d} E \bigl[ \bigl(Z^\epsilon_{nm}
\bigr)^{2} \bigr] \\
&&\qquad=\int g(x_{1}, \ldots, x_{K-1},
y_{0}-\epsilon v, y_{1}, \ldots, y_{L-1})
\mu_{2}(y_{0}, y_{1}, \ldots, y_{L-1})
\\
&&\hspace*{8pt}\qquad\quad{}\times K(v)^{2} \prod_{i=1}^{K-1}
p(s_{i-1}, x_{i-1}, s_{i}, x_{i})
p(s_{K-1}, x_{K-1}, s_{K}, y_{0} -
\epsilon v)
\\
&&\hspace*{8pt}\qquad\quad{}\times\prod_{i=1}^{L} q(
\widehat{t}_{i-1}, y_{i}, \widehat{t}_{i},
y_{i-1}) \times dx_{1} \cdots dx_{K-1} \,dv
\,dy_{0} \,dy_{1} \cdots dy_{L-1}.
\end{eqnarray*}
For $\epsilon\to0$ the right-hand side gives the statement from the lemma.

For the difference, consider
\begin{eqnarray*}
r_{\epsilon}^{(1,1)} &\coloneqq&\int K(v)^{2} \bigl[
g(x_{1},\ldots, x_{K-1}, y_{0}-\epsilon v,
y_{1}, \ldots, y_{L-1})^{2} \\
&&\hspace*{19pt}\qquad {}\times p(s_{K-1},
x_{K-1}, s_{K}, y_{0}-\epsilon v)
\\
&&\hspace*{8pt}\qquad\quad{}- g(x_{1},\ldots, x_{K-1}, y_{0},
y_{1}, \ldots, y_{L-1})^{2}\\
&&\hspace*{77pt}\qquad\quad{}\times p(s_{K-1},
x_{K-1}, s_{K}, y_{0}) \bigr] \,dv.
\end{eqnarray*}
Following the procedure established in the previous lemmas, we obtain
\[
\bigl| r^{(1,1)}_{\epsilon}\bigr| \le C s_{\epsilon}^{(1,1)}(x_{K-1},
y_{0}),
\]
and by the argument used in the proof of Lemma \ref{lem:Zn-nprime}, we obtain
transition densities function $\widetilde{p}(s_{0},x,t^{\ast},y_{0})$ and
$\widetilde{q}(t^{\ast},y_{0}, T,y)$ such that
%
%e35 #&#
\begin{eqnarray}
\label{eq:Znm}&&\Bigl\llvert \epsilon^{d} \mathbb{E}\bigl[
\bigl(Z^\epsilon_{nm}\bigr)^{2}\bigr] - \lim
_{\epsilon\to0} \epsilon^{d} \mathbb{E} \bigl[
\bigl(Z_{nm}^\epsilon\bigr)^{2} \bigr] \Bigr\rrvert
\nonumber
\\[-8pt]
\\[-8pt]
\nonumber
&&\qquad
\le C \epsilon^{2} \int\widetilde{p}\bigl(s_{0},x,t^{\ast},y_{0}
\bigr) \widetilde{q}\bigl(t^{\ast},y_{0}, T,y\bigr)
\,dy_{0}. %\qedhere
\end{eqnarray}
\upqed\end{pf}

In what follows, we simplify the notation by the following conventions:
\begin{itemize}
\item the constant in Theorem \ref{bias} is denoted by $C_0$, that is,
$|h_{\epsilon}- h| \le C_0 \epsilon^{2}$;
\item for $m \neq m^{\prime}$, we set $\mathbb{E} [ Z^\epsilon_{nm}
Z^\epsilon_{nm^{\prime}} ] \eqqcolon h_{\epsilon}^{(1,2)}$ and denote
the constant for the difference by $C_{1,2}$, that is, $\llvert
h^{(1,2)}_{\epsilon}- h^{(1,2)}_{0} \rrvert  \le C_{1,2} \epsilon^{2}$;
\item for $n \neq n^{\prime}$, we set $\mathbb{E} [ Z^\epsilon_{nm}
Z^\epsilon_{n^{\prime}m} ] \eqqcolon h_{\epsilon}^{(2,1)}$ and denote
the constant for the difference by $C_{2,1}$, that is, $\llvert
h^{(2,1)}_{\epsilon}- h^{(2,1)}_{0} \rrvert  \le C_{2,1} \epsilon^{2}$;
\item we set $\epsilon^{d} \mathbb{E} [ (Z^\epsilon_{nm})^{2} ]
\eqqcolon h_{\epsilon}^{(1,1)}$ and denote the constant for the difference
by $C_{1,1}$, that is, $\llvert  h^{(1,1)}_{\epsilon}- h^{(1,1)}_{0}
\rrvert  \le
C_{1,1} \epsilon^{2}$.
\end{itemize}

%le4.13 #&#
\begin{lemma}
\label{lem:variance} The variance of the estimator is given by
\[
\operatorname{Var} \widehat{h}_{\epsilon,M,N} = \frac{1-M-N}{NM}
h_{\epsilon
}^{2} + \frac{M-1}{NM} h_{\epsilon}^{(1,2)}
+ \frac{N-1}{NM} h_{\epsilon
}^{(2,1)} + \frac{\epsilon^{-d}}{NM}
h_{\epsilon}^{(1,1)}.
\]
\end{lemma}
%
%re4.14 #&#
\begin{remark}
Lemma \ref{lem:variance} gives a clarification of the intuitive fact that
the variance of $\widehat{h}_{\epsilon,M,N}$ explodes as $\epsilon\to0$
(and, hence, $K_\epsilon\to\delta_0$). Indeed, as all the
$h_\epsilon^{(\cdot)}$ terms have a finite limit, the explosion is exclusively
caused by the contribution of $\mathbb{E}[(Z^\epsilon_{nm})^2] =
\epsilon^{-d} h^{(1,1)}_\epsilon$. Finally, the exploding term
$\epsilon^{-d}$ will be compensated by the factor $1/(NM)$.
\end{remark}
\begin{pf*}{Proof of Lemma \ref{lem:variance}}
The result follows immediately by \eqref{eq:def-Z}, independence of
$Z^\epsilon_{nm}$ and $Z^\epsilon_{n^{\prime}m^{\prime}}$ when both
$n\neq
n^{\prime}$ and $m \neq m^{\prime}$ and the notation introduced above, noting
that $E[Z^\epsilon_{nm}] = h_{\epsilon}$.
\end{pf*}

We immediately obtain the following:\vadjust{\goodbreak}
%
%le4.15 #&#
\begin{lemma}
\label{lem:mse-h}
We assume Conditions \ref{ass:bound-density},  \ref{ass:kernel-order} and
\ref{ass:convenience} hold. Then the
mean square error of the estimator $\widehat{h}_{\epsilon,M,N}$ introduced
in \eqref{eq:hat-heps} for the term $h$ defined in~\eqref{eq:def-h}
satisfies
\begin{eqnarray*}
&&\mathbb{E} \bigl[ ( \widehat{h}_{\epsilon,M,N} - h ) ^{2} \bigr]\\
&&\qquad \le
\frac{1-N-M}{NM} h^{2} + \frac{M-1}{NM} h^{(1,2)}_{0}
+ \frac
{N-1}{NM} h^{(2,1)}_{0} + \frac{\epsilon^{-d}}{NM}
h^{(1,1)}_{0}
\\
&&\qquad\quad{}+ \frac{\epsilon^{-d+2}}{NM} C_{1,1} + \epsilon^{2} \biggl[ 2
\frac{1-N-M}{NM} C h + \frac{M-1}{NM} C_{1,2} + \frac{N-1}{NM}
C_{2,1} \biggr]\\
&&\quad\qquad{} + \frac{(N-1)(M-1)}{NM} C_0^{2}
\epsilon^{4}.
\end{eqnarray*}
\end{lemma}

Similar to \citet{MSS1}, we can now choose $N=M$ and the bandwidth
$\epsilon$
so as to obtain convergence proportional to $N^{-1/2}$ in RMSE-sense.

%th4.16 #&#
\begin{theorem}
\label{thr:mse-h-order} Assume Conditions \ref{ass:bound-density},
 \ref{ass:kernel-order} and \ref{ass:convenience} and set $M = N$, and
$\epsilon= \epsilon_{N}$ dependent on $N$.
\begin{itemize}
\item If $d \le4$, choose $\epsilon_{N} = C N^{-\alpha}$ for some $1/4
\le\alpha\le
1/d$. Then we have $\mathbb{E} [  ( \widehat{h}_{\epsilon
_{N},N,N} -
h  ) ^{2}  ] = \mathcal{O}(N^{-1})$, so we achieve the optimal
convergence rate $1/2$.
\item For $d > 4$, choose $\epsilon_{N} = C N^{-2/(4+d)}$, and we obtain
$\mathbb{E} [  ( \widehat{h}_{\epsilon_{N},N,N} - h  ) ^{2}
] = \mathcal{O}(N^{-8/(4+d)})$.
\end{itemize}
\end{theorem}
\begin{pf}
Insert $M=N$ and the respective choice of $\epsilon_{N}$ in
Lemma \ref{lem:mse-h}.
\end{pf}

%re4.17 #&#
\begin{remark}
\label{rem:higher-order-kernel}
By replacing the kernel $K$ by \emph{higher order} kernels,\footnote{Recall that the order of a kernel $K$ is the order of the lowest order
(nonconstant) monomial $f$ such that $\int f(v) K(v) \,dv \neq0$.} one
could retain the convergence rate $1/2$ even in higher dimensions, as higher
order kernels lead to higher order estimates (in $\epsilon$) in
Lemmas \ref{lem:Z-m-mprime}, \ref{lem:Zn-nprime} and~\ref{lem:Znm}.
\end{remark}

So far, we have only computed the quantity $h$ as given in \eqref{eq:def-h}.
However, finally we want to compute the conditional expectation
\[
H \coloneqq\mathbb{E} \bigl[  g \bigl( X_{s_{0},x}(s_{1}),
\ldots, X_{s_{0},x}(t_{L-1}) \bigr) \rrvert X_{s_{0},x}(T)=y
\bigr].
\]
As $H = \frac{h}{p(s_{0},x,T,y)}$ with $h$ defined in \eqref{eq:def-h}, we
need to divide the estimator for $h$ by an appropriate estimator for
$p(s_{0},x,T,y)$---in fact, we choose the forward reverse estimator
with $g
\equiv1$. Note that we have assumed that $p(s_{0}, x, T, y) > 0$. To
rule out large error contributions when the denominator is small, we will
discard experiments which give too small estimates for the transition density.
More precisely, we choose our final estimator to be
%
%e36 #&#
%e37 #&#
\begin{eqnarray}
\label{eq:H-hat-def}
\widehat{H}_{\epsilon,M,N} &\coloneqq&\frac{\sum_{n=1}^{N} \sum_{m=1}^{M}
g ( X_{s_{1}}^{n}, \ldots, X_{s_{K}}^{n}, Y_{\widehat{t}_{L-1}}^{m},
\ldots, Y_{\widehat{t}_{1}}^{m}  ) K ( {(Y_{\widehat{t}_L}^{m} -
X_{t^{\ast}}^{n})}/{\epsilon}  )
\mathcal{Y}_{\widehat{t}_L}^{m}}{\sum_{n=1}^{N} \sum_{m=1}^{M} K (
{(Y_{\widehat{t}_L}^{m} - X_{t^{\ast}}^{n})}/{\epsilon}  )
\mathcal{Y}_{\widehat{t}_L}^{m}}
\nonumber
\\[-8pt]
\\[-8pt]
\nonumber
&&{}\times\mathbf{1}_{({1}/{(NM)}) \epsilon^{-d} \sum_{n=1}^{N}
\sum_{m=1}^{M}
K (( {Y_{\widehat{t}_L}^{m} - X_{t^{\ast}}^{n})}/{\epsilon}
)
\mathcal{Y}_{\widehat{t}_L}^{m} > \overline{p}/2},
\end{eqnarray}
where $\overline{p} > 0$ is a lower bound for $p(s_{0},x,T,y)$ (for fixed
$s_{0},x,T,y$), which is assumed to be known.\footnote{In practice,
such a
lower bound could be achieved by running an independent estimation for
$p(s_{0},x,T,y)$ and then taking a value at the lower end of a required
confidence interval. See Remark~\ref{rem:nullfolge} below for a different
version of the theorem. In any case, our numerical experiments suggest that
the cut-off can be safely omitted in practice. Keep in mind, however, that
the ratio of the asymptotic distributions for numerator and
denominator may
not have finite moments.}

%th4.18 #&#
\begin{theorem}
\label{thr:mse-H} Assume Conditions \ref{ass:bound-density},
\ref{ass:kernel-order} and \ref{ass:convenience} and set $M = N$ and
$\epsilon= \epsilon_{N}$ dependent on $N$.
\begin{itemize}
\item If $d \le4$ (or $d>4$ and higher order kernels are used), choose
$\epsilon_{N} = C N^{-\alpha}$, $1/4 \le\alpha\le1/d$. Then we have
$\mathbb{E} [  ( \widehat{H}_{\epsilon_{N},N,N} - H  ) ^{2}
] = \mathcal{O}(N^{-1})$, so we achieve the optimal convergence rate
$1/2$.
\item For $d > 4$, choose $\epsilon_{N} = C N^{-2/(4+d)}$, and we obtain
$\mathbb{E} [  ( \widehat{H}_{\epsilon_{N},N,N} - H  ) ^{2}
] = \mathcal{O}(N^{-8/(4+d)})$.
\end{itemize}
\end{theorem}

\begin{pf}
Let $X_{N} \coloneqq\widehat{h}_{\epsilon_{N},N,N}$, and, similarly, let
\[
Y_{N} \coloneqq\frac{1}{N^{2}} \epsilon_{N}^{-d}
\sum_{n=1}^{N} \sum
_{m=1}^{N} K \biggl( \frac{Y_{\widehat{t}_L}^{m} - X_{t^{\ast}}^{n}}{\epsilon_{N}} \biggr)
\mathcal{Y}_{\widehat{t}_L}^{m}
\]
denote the estimator in the denominator, including the normalization factor.
Moreover, let $X \coloneqq h$ as defined in \eqref{eq:def-h} and let $Y
\coloneqq p(s_{0}, x, T, y)$. Then we have already established in
Theorem \ref{thr:mse-h-order}
that
\begin{eqnarray*}
\mathbb{E} \bigl[ |X_{N} - X|^{2} \bigr] & =& \mathcal{O}
\bigl(N^{-p}\bigr),
\\
\mathbb{E} \bigl[ |Y_{N} - Y|^{2} \bigr] & =& \mathcal{O}
\bigl(N^{-p}\bigr),
\end{eqnarray*}
where $p = 1$ for $d \le4$ and $p = \frac{8}{4+d}$ when $d > 4$.
Moreover, we
have obtained in Lemma \ref{lem:variance} that $\operatorname{Var}
X_{N} =
\mathcal{O}(N^{-p})$ and $\operatorname{Var} Y_{N} = \mathcal{O}(N^{-p})$.

We will now estimate the mean square error for the quotient by
splitting it
into two contributions, depending\vadjust{\goodbreak} on whether $Y_{N}$ is small or large. To
this end, let
\[
\zeta_{N} \coloneqq\frac{X_{N}}{Y_{N}} \mathbf{1}_{Y_{N} > D_{N}}%
\]
for a constant $D_{N}$ to be specified below satisfying $D_{N} <
\mathbb{E}[Y_{N}]$ (in fact, for $N$ large enough, this constant may be
chosen to be $\overline{p}/2$). Then we have
%
%e38 #&#
\begin{eqnarray}
\label{eq:2}%
\mathbb{E} \biggl[ \biggl( \frac{X_{N}}{Y_{N}} -
\frac{X}{Y} \biggr) ^{2} \mathbf{1}_{Y_{N} > D_{N}} \biggr] & =&
\mathbb{E} \biggl[ \frac{ (
X_{N} Y - Y_{N} X  ) ^{2}}{(Y_{N}Y)^{2}} \mathbf{1}_{Y_{N} > D_{N}} \biggr]
\nonumber
\\
& \le&\frac{\mathbb{E} [  ( Y(X_{N}-X) + X(Y-Y_{N})  ) ^{2}
] }{Y^{2} D_{N}^{2}}
\nonumber
\\[-8pt]
\\[-8pt]
\nonumber
& \le&2 \frac{Y^{2} \mathbb{E}[(X_{N}-X)^{2}] + X^{2} \mathbb{E} [
(Y-Y_{N})^{2}  ] }{Y^{2} D_{N}^{2}}
\\
& \le&\frac{C^{1}_{X,Y}}{D_{N}^{2} N^{p}},\nonumber
\end{eqnarray}
where we used the estimates on the MSEs for numerator and denominator.
On the
other hand, we have, using that $D_{N} < \mathbb{E} Y_{N}$, Chebyshev's
inequality and our estimate on the variance of $Y_{N}$,
%
%e39 #&#
\begin{eqnarray}
\label{eq:4} \mathbb{P} ( Y_{N} \le D_{N} ) & =&
\mathbb{P} ( Y_{N} - \mathbb{E} Y_{N} \le D_{N} -
\mathbb{E} Y_{N}; Y_{N} \le\mathbb{E} Y_{N} )
\nonumber
\\
& \le&\mathbb{P}\bigl ( |Y_{N} - \mathbb{E} Y_{N}| \ge
\mathbb{E}Y_{N} - D_{N} \bigr)
\nonumber
\\[-8pt]
\\[-8pt]
\nonumber
& \le&\frac{\operatorname{Var} Y_{N}}{ ( \mathbb{E} Y_{N} - D_{N}
) ^{2}}
\\
& \le&\frac{C^{2}_{Y}}{ ( \mathbb{E} Y_{N} - D_{N}  ) ^{2} N^{p}}.\nonumber
\end{eqnarray}
Finally, consider
%
%e40 #&#
\begin{eqnarray}
\label{eq:6} \mathbb{E} \biggl[ \biggl( \zeta_{N} -
\frac{X}{Y} \biggr) ^{2} \biggr] & = &\mathbb{E} \biggl[ \biggl(
\zeta_{N} - \frac{X}{Y} \mathbf{1}_{Y_{N} >
D_{N}} -
\frac{X}{Y} \mathbf{1}_{Y_{N} \le D_{N}} \biggr) ^{2} \biggr]
\nonumber
\\
& =& \mathbb{E} \biggl[ \biggl( \zeta_{N} - \frac{X}{Y}
\mathbf{1}_{Y_{N}>D_{N}} \biggr)^{2} \biggr] + \frac{X^{2}}{Y^{2}}
\mathbb{P} ( Y_{N} \le D_{N} )
\\
& \le&\frac{C^{1}_{X,Y}}{D_{N}^{2} N^{p}} + \frac{C^{2}_{Y} X^{2}}{ (
\mathbb{E} Y_{N} - D_{N}  )^{2} Y^{2} N^{p}},\nonumber
\end{eqnarray}
where we have combined \eqref{eq:2} and \eqref{eq:4}. Now choose $D_{N}
\coloneqq\overline{p}/2$ for $N$ large enough. As $\mathbb{E} Y_{N}
\mathop{\xrightarrow}\limits^{N\to\infty} Y$, \eqref{eq:6} implies that
%
%e41 #&#
\begin{equation}
\label{eq:7}\mathbb{E} \biggl[ \biggl( \zeta_{N} -
\frac{X}{Y} \biggr) ^{2} \biggr] = \mathcal{O}
\bigl(N^{-p}\bigr). %\qedhere
\end{equation}
\upqed\end{pf}

%re4.19 #&#
\begin{remark}
\label{rem:nullfolge} Alternatively, we could replace the cut-off
$\overline{p}/2$ in \eqref{eq:H-hat-def} by some sequence $D_{N}
\mathop{\xrightarrow}\limits^{N\to\infty} 0$. In that case, the MSE of the estimator is of
order $\mathcal{O}(N^{-p}/D_{N}^{2})$, which can be chosen as close to
$\mathcal{O}(N^{-p})$ as desired by proper choices of (slowly convergent)
sequences $D_{N}$. Note that finally $\mathbb{E} Y_{N} > D_{N}$ in the proof
of Theorem \ref{thr:mse-H}, as $\mathbb{E}Y_{N} \mathop{\xrightarrow}\limits^{N\to\infty}
p(s_{0}, x, T, y) > 0$ by assumption.
\end{remark}

%s4.2 #&#
\subsection{Forward-reverse estimators for conditioning on a set}
\label{sec:forw-reverse-estim-set}

In Theorem \ref{thr:cond-dist-set} and Corollary \ref
{cor:cond-dist-comp} we
have derived a representation of the conditional expectation of a functional
$g$ of the process $X$ given that $X_{T} \in A$ (for a Borel set $A$ with
positive probability) or given $X_{T}^{1} = y^{1}, \ldots, X_T^{d'} =
y^{d'}$. In analogy to the first part of this section, one can construct
Monte Carlo estimators for these conditional expectations and analyze their
bias and variance. In what follows, we assume that $A$ is either a general
Borel set with positive probability or an affine surface, that is, we
treat both
cases distinguished above together.

Recall that we represented the conditional expectation as
\begin{eqnarray*}
&&\lim_{\epsilon\downarrow0} \mathbb{E} \biggl[ g \bigl( X_{s_{0},x}(s_{1}),
\ldots, X_{s_{0},x}(s_{K}), Y_{\xi;T}\bigl(
\widehat{t}_{L-1}, \ldots Y_{t^{\ast
},\xi}(\widehat{t}_{1})
\bigr)\bigr)
\\
&&\hspace*{59pt}\quad{}\times\epsilon^{-d} K \biggl( \frac{Y_{\xi;T}(\widehat{t}_L) - X_{s_{0},x}(t^{\ast})}{\epsilon} \biggr)
\frac{\mathcal{Y}_{\xi;T}(\widehat{t}_L)}{\varphi(\xi)} \biggr]
\\
&&\qquad= \int_{A} p(s_{0},x,T,y)
\lambda_{A}(dy) \mathbb{E} \bigl[  g \bigl(
X_{s_{0},x}(s_{1}), \ldots, X_{s_{0},x}(t_{L-1})
\bigr) \rrvert X_{s_{0},x}(T) \in A \bigr],
\end{eqnarray*}
where $\xi$ is an independent random variable taking values in $A$ with
density $\varphi$ with respect to $\lambda_{A}$. In order to arrive at an
estimator with bounded variance, we need to restrict the choice of
$\varphi$
and, consequently, $\xi$.

%co4.20 #&#
\begin{condition}
\label{ass:varphi} The density $\varphi$ has (strictly) super-Gaussian tails,
that is, there are constants $C, \gamma, \delta> 0$ such that
\[
\varphi(v)^{-1} \le C \exp \bigl( \gamma|v|^{2-\delta} \bigr),\qquad v
\in A.
\]
\end{condition}

We define the following Monte Carlo estimator for the conditional
expectation
%
%e42 #&#
%e43 #&#
\fontsize{10pt}{\baselineskip}\selectfont
\makeatletter
\def\tagform@#1{\normalsize\maketag@@@{(\ignorespaces#1\unskip\@@italiccorr)}}
\makeatother
\begin{eqnarray}
\label{eq:H-comp-def}&&\widehat{H}_{\epsilon,M,N}^{\xi}\nonumber\\
&&\qquad\coloneqq
\frac
{\sum_{n=1}^{N}
\sum_{m=1}^{M} g ( X_{s_{1}}^{n}, \ldots, X_{s_{K}}^{n}, Y_{\widehat
{t}%
_{L-1}}^{m}, \ldots, Y_{\widehat{t}_{1}}^{m}  ) K (
{(Y_{\widehat{t}_L}^{m} -
X_{t^{\ast}}^{n})}/{\epsilon}  ) ({\mathcal{Y}_{\widehat
{t}_L}^{m}}/{\varphi
(\xi^{m})})}{\sum_{n=1}^{N} \sum_{m=1}^{M} K ( {(Y_{\widehat
{t}_L}^{m} -
X_{t^{\ast}}^{n})}/{\epsilon}  ) ({\mathcal{Y}_{\widehat
{t}_L}^{m}}/{\varphi
(\xi^{m})})}
\\
&&\qquad\quad{}\times\mathbf{1}_{({1}/{(NM)} )\epsilon^{-d} \sum_{n=1}^{N} \sum_{m=1}^{M}
K ( {(Y_{\widehat{t}_L}^{m} - X_{t^{\ast}}^{n})}/{\epsilon}
) ({\mathcal{Y}_{\widehat{t}_L}^{m}}/{\varphi(\xi^{m})}) > \overline{p}/2},\nonumber
\end{eqnarray}
\normalsize
\hspace*{-3pt}where $(X_{s_{1}}^{n}, \ldots, X_{s_{K}}^{n})$, $1 \le n \le N$, are
independent samples from the solution of the forward process $X$
started at
$X_{s_{0}} = x$ and $(Y^{m}_{\widehat{t}_{L}}, \ldots,
Y^{m}_{\widehat{t}_{1}})$ together with $\mathcal{Y}^{m}_{\widehat
{t}_L}$, $1
\le m \le M$, are independent samples from the reverse process $(Y,
\mathcal{Y})$ started at $Y^{m}_{0} = \xi^{m}$, $\mathcal{Y}^{m}_{0} =
1$, for
an independent sequence of samples $\xi^{m}$ from the distribution~$\xi$.
Apart from the term $\varphi(\xi^{m})$, the difference to
estimator \eqref{eq:H-hat-def} is the randomness of the initial values
of the
reverse process. Again, $p(s_{0},x,T,y) > \overline{p} > 0$, and
Remark~\ref{rem:nullfolge} applies. The analysis of \eqref{eq:H-comp-def},
however, works along the lines of the analysis
of \eqref{eq:H-hat-def}. Indeed, in all the expectations considered in
Theorem \ref{bias} and in Lemmas~\ref{lem:Z-m-mprime}--\ref{lem:Znm}, we obtain
the same kind of results by the following steps:

\begin{longlist}[(1)]
\item[(1)] condition on $\xi$ and pull out the factor $\varphi(\xi
)^{-1}$ (possibly
with indices $m$ and/or $m^{\prime}$);
\item[(2)] use the results obtained in Section~\ref{sec:forw-reverse-estim-state},
with constants depending on the value of $\xi$;
\item[(3)] move $\varphi(\xi)^{-1}$ back in and take the expectation in
$\xi$.
\end{longlist}

%th4.21 #&#
\begin{theorem}
\label{thr:mse-comp} Set $M=N$ and assume Condition \ref{ass:varphi}
and, as
usual, Condition \ref{ass:bound-density}, \ref{ass:kernel-order}
and \ref{ass:convenience}.

\begin{itemize}
\item If $d \le4$, choose $\epsilon_{N} = C N^{-\alpha}$, $1/4 \le
\alpha\le
1/d$. Then the MSE of the forward-reverse estimator
$\widehat{H}_{\epsilon,M,N}^{\xi}$ is $\mathcal{O}%
(N^{-1})$.

\item For $d > 4$, choose $\epsilon_{N} = C N^{-2/(4+d)}$. Then the MSE
of the
forward-reverse estimator $\widehat{H}_{\epsilon,M,N}^{\xi}$ is
$\mathcal{O}%
(N^{-8/(4+d)})$.
\end{itemize}
\end{theorem}

\begin{pf}
In this proof, the constant $C$ may change from line to line. Define
\begin{eqnarray*}
h^{\xi} & \coloneqq &\int_{A}p(s_{0},x,T,y)
\lambda_{A}(dy)\cdot\mathbb{E} \bigl[  g \bigl(
X_{s_{0},x}(s_{1}),\ldots,X_{s_{0},x}(t_{L-1})
\bigr) \rrvert X_{s_{0},x}(T)\in A \bigr],
\\
h^{\xi}_{\epsilon} & \coloneqq&\mathbb{E} \biggl[g \bigl(
X_{s_{0},x}(s_{1} ),\ldots,X_{s_{0},x}
\bigl(t^{\ast}\bigr),Y_{\xi;T}(\widehat{t}_{L-1}),\ldots
,Y_{\xi;T}(\widehat{t}_{i}) \bigr)
\\
&&\hspace*{80pt}{} \times K_{\epsilon} \bigl( Y_{\xi;T}(\widehat{t}_L)-X_{s_{0},x}
\bigl(t^{\ast
}\bigr) \bigr) \frac{\mathcal{Y}_{\xi;T}(\widehat{t}_L)}{\varphi(\xi)} \biggr],
\\
Z^{\epsilon,\xi}_{nm} & \coloneqq&\frac{1}{\epsilon^{d}} g \bigl(
X_{s_{0},x}^{n}(s_{1}), \ldots, X_{s_{0},x}^{n}(s_{K}),
Y_{\xi^{m};T}^{m}(\widehat{t}_{L-1}), \ldots,
Y_{\xi^{m};T}^{m}(\widehat{t}_{1}) \bigr)
\\
&&{} \times K \biggl( \frac{Y_{\xi^{m};T}^{m}(\widehat{t}_L) - X_{s_{0},x}^{n}(t^{\ast
})}{\epsilon} \biggr) \frac{\mathcal{Y}^{m}_{\xi^{m};T}(\widehat{t}_L)}{\varphi(\xi^{m})},
\end{eqnarray*}
and notice that the result will follow if we can establish the bounds of
Theorem \ref{bias} and Lemmas \ref{lem:Z-m-mprime}, \ref{lem:Zn-nprime}
and \ref{lem:Znm} for $h$, $h_{\epsilon}$ and $Z^\epsilon_{nm}$
replaced by
$h^{\xi}$, $h^{\xi}_{\epsilon}$ and $Z^{\epsilon,\xi}_{nm}$, respectively.

For the bias, \eqref{eq:bias-bound-explicit} implies a bound $|h(y) -
h_{\epsilon}(y)| \le C \epsilon^{2} \widetilde{p}(s_{0},x,T,y)$ for some
density $\widetilde{p}$ in $y$, where we make the dependence of $h$ and
$h_{\epsilon}$ on $y$ explicit. Consequently, conditioning on $\xi$
first, we
have
%
%e44 #&#
\begin{eqnarray}
\label{eq:bias-bound-xi}%
\bigl|h^{\xi}- h^{\xi}_{\epsilon}\bigr| &
=& \biggl\llvert \mathbb{E} \biggl[ \frac
{h(\xi)-h_{\epsilon}(\xi)}{\varphi(\xi)} \biggr] \biggr\rrvert
\nonumber
\\
& \le&\mathbb{E} \biggl[ \frac{|h(\xi) - h_{\epsilon}(\xi)|}{\varphi(\xi)} \biggr]
\nonumber
\\[-8pt]
\\[-8pt]
\nonumber
& \le& C \epsilon^{2} \int\frac{\widetilde{p}(s_{0},x,T,\xi)}{\varphi(\xi)} \varphi(\xi) \,d\xi
\\
& \le& C \epsilon^{2}.\nonumber
\end{eqnarray}

Similarly, using the estimate from Lemma \ref{lem:Z-m-mprime},
denoting $z^\epsilon_{n,m,m^{\prime}}(y,y^{\prime}) \coloneqq\mathbb
{E} [
Z^\epsilon_{nm} Z^\epsilon_{nm^{\prime}} ] $, where we assume
$Y^{m} =
Y^{m}_{y;T}$ and $Y^{m^{\prime}} =
Y^{m^{\prime}}_{t^{\ast},y^{\prime}}$, we get, using a simple adaptation
of \eqref{eq:Zm-mprime} for different terminal values $y$ and $y^{\prime}$,
%
%e45 #&#
\begin{eqnarray}
\label{eq:Zm-mprime-bound-xi} &&\bigl| \mathbb{E} \bigl[ Z_{nm}^{\epsilon,\xi}
Z_{nm^{\prime}}^{\epsilon,\xi} \bigr] -  \mathbb{E} \bigl[
Z_{nm}^{\epsilon,\xi}Z_{nm^{\prime}}^{\xi} \bigr] |
_{\epsilon=0} \bigr\rrvert\nonumber \\
&&\qquad\le\mathbb{E} \biggl[ \frac{\llvert  z^\epsilon
_{n,m,m^{\prime
}}(\xi^{m},\xi^{m^{\prime}}) -
z^\epsilon_{n,m,m^{\prime}}(\xi^{m},\xi^{m^{\prime}}) \rrvert  _{\epsilon
=0} \rrvert
}{\varphi(\xi^{m}) \varphi(\xi^{m^{\prime}})} \biggr]
\nonumber
\\
&&\qquad \le C \epsilon^{2} \mathbb{E} \biggl[ \frac{\int\widetilde{p}(s_{0}, x,
t^{\ast},x_{K}) \widetilde{q}(t^{\ast},x_{K}, T, \xi^{m}) \widetilde
{q}(t^{\ast},x_{k},T,\xi^{m^{\prime}}) \,dx_{K}}{\varphi(\xi^{m}) \varphi
(\xi^{m^{\prime}})} \biggr]
\nonumber
\\[-8pt]
\\[-8pt]
\nonumber
&&\qquad = C\epsilon^{2} \int\widetilde{p}\bigl(s_{0},x,t^{\ast},x_{K}
\bigr) \widetilde {q}\bigl(t^{\ast},x_{K}, T, y\bigr)\\
&&\hspace*{28pt}\qquad\quad{}\times
\widetilde{q}\bigl(t^{\ast},x_{k},T, y^{\prime}\bigr)
\,dx_{K} \lambda_{A}(dy) \lambda_{A}
\bigl(dy^{\prime}\bigr)
\nonumber
\\
&&\qquad \le C \epsilon^{2}.\nonumber
\end{eqnarray}

Adopting the above notation for the case $n\neq n^{\prime}$ covered in
Lemma \ref{lem:Zn-nprime} and using \eqref{eq:Zn-nprime}, we get
\begin{eqnarray*}
&&\bigl\llvert \mathbb{E} \bigl[ Z_{nm}^{\epsilon,\xi}
Z_{n^{\prime}m}^{\epsilon,\xi} \bigr] -  \mathbb{E} \bigl[
Z_{nm}^{\epsilon,\xi}Z_{n^{\prime}m}^{\epsilon,\xi
} \bigr]
\rrvert _{\epsilon=0}\bigr \rrvert\\
&&\qquad \le\mathbb{E} \biggl[ \frac{\llvert  z^\epsilon
_{n,n^{\prime
},m}(\xi^{m},\xi^{m}) -  z^\epsilon_{n,n^{\prime},m}(\xi^{m},\xi
^{m}) \rrvert
_{\epsilon=0} \rrvert  }{\varphi(\xi^{m}) \varphi(\xi^{m})}
\biggr]
\\
&&\qquad \le C \epsilon^{2} \int\frac{\widetilde{p}(s_{0},x,t^{\ast},y_{1})
\widetilde{q}(t^{\ast},y_{1}, T, y)}{\varphi(y)} \,dy_{1}
\lambda_{A}(dy).
\end{eqnarray*}
By assumption the density $\int\widetilde{p}(s_{0},x,t^{\ast},y_{1})
\widetilde{q}(t^{\ast},y_{1}, T, y) \,dy_{1}$ has Gaussian tails, whereas
$\varphi$ was assumed to have strictly sub-Gaussian tails. This implies that
the above integral is finite, and we get the bound
%
%e46 #&#
\begin{equation}
\label{eq:Zn-nprime-bound-xi}\bigl\llvert \mathbb{E} \bigl[ Z_{nm}^{\epsilon
,\xi
}Z_{n^{\prime}m}^{\epsilon,\xi}
\bigr] -  \mathbb{E} \bigl[ Z_{nm}^{\epsilon,\xi
}Z_{n^{\prime}m}^{\epsilon,\xi}
\bigr] \rrvert _{\epsilon=0} \bigr\rrvert \le C \epsilon^{2}.
\end{equation}

In a similar way, using \eqref{eq:Znm}, we get the bound
%
%e47 #&#
\begin{equation}
\label{eq:Znm-bound-xi}\Bigl\llvert \epsilon^{d} \mathbb{E} \bigl[
\bigl(Z_{nm}^{\epsilon,\xi
}\bigr)^{2} \bigr] - \lim
_{\epsilon\to0} \epsilon^{d} \mathbb{E} \bigl[
\bigl(Z_{nm}^{\epsilon,\xi}\bigr)^{2} \bigr] \Bigr\rrvert
\le C \epsilon^{2}.
\end{equation}

The respective versions of Lemmas \ref{lem:variance}, \ref
{lem:mse-h} and
Theorem \ref{thr:mse-h-order} follow immediately from the
bounds \eqref{eq:bias-bound-xi}, \eqref{eq:Zm-mprime-bound-xi},
\eqref{eq:Zn-nprime-bound-xi} and \eqref{eq:Znm-bound-xi}, and we can repeat
the proof of Theorem \ref{thr:mse-H}, arriving at the conclusion.
\end{pf}

We again stress that the nonoptimal complexity rate in
Theorem \ref{thr:mse-comp} can be improved to the optimal one even for
$d > 4$
by Remark \ref{rem:higher-order-kernel}.

%s4.3 #&#
\subsection{Limitations of the forward-reverse estimator}
\label{sec:limit-forw-reverse}

Theorems \ref{thr:mse-H} and \ref{thr:mse-comp} above present the asymptotic
analysis of the MSE for the forward-reverse estimator. In practice, for many
methods with very good asymptotic rates, limitations arise due to potentially
high constants, and the forward-reverse estimator is no exception. In fact,
this can be already seen in a very simple example, where all the
estimates can
be given explicitly.

For $s_0 = 0 < t^\ast< T$, consider the one-dimensional Ornstein--Uhlenbeck
process
%
%e48 #&#
\begin{equation}
\label{eq:ou-example} dX_{0,x}(t) = -\alpha X_{0,x}(t) \,dt +
dB_t
\end{equation}
for $\alpha>0$. The corresponding reverse process satisfies
%
%e49 #&#
\begin{equation}
\label{eq:ou-example-reverse} dY_{y;T}(t) = \alpha Y_{y;T}(t) \,dt +
dW_t
\end{equation}
for a Brownian motion $W_t$. Moreover, $\mathcal{Y}_{y;T}(T-t^\ast) =
e^{\alpha(T-t^\ast)}$. We first discuss the estimator $\widehat
{h}_{\epsilon,N,N}$
introduced in \eqref{eq:hat-heps} for the numerator of the forward-reverse
estimator $\widehat{H}_{\epsilon,N,N}$ for $g \equiv1$ with $K = L =
1$. Of
course, we expect that the findings for this special case carry over to
situations with nonconstant $g$ and $K,L \ge1$.

After elementary but tedious calculations
[Milstein, Schoenmakers and
Spokoiny (\citeyear{MSS1}), Section~4] one
arrives at
%
%e50 #&#
\begin{equation}
\label{eq:ou-example-mean}\qquad  \mathbb{E} [ \widehat{h}_{\epsilon,N,N} ] = \f{1} {\sqrt
{2 \pi \bigl( \epsilon^2 e^{-2\alpha(T-t^\ast)} +
\sigma^2_T \bigr)}} \exp \biggl( - \f{ (
e^{-\alpha T}x - y )^2} {2 ( \epsilon^2
e^{-2\alpha(T-t^\ast)} + \sigma^2_T )} \biggr)
\end{equation}
and
%
%e51 #&#
%e52 #&#
%e53 #&#
%e54 #&#
\begin{eqnarray}
\label{eq:ou-example-var} \operatorname{Var} \widehat{h}_{\epsilon,N,N}& = &-
 \f{2N-1} {2\pi
N^2 (B+\sigma^2_T)} \exp \biggl( -
\f{A} {B+\sigma^2_T} \biggr)\nonumber
\\
&&{}+ \f{N-1} {2\pi N^2 \sqrt{B+\sigma^2_{T-t^\ast}}
\sqrt{B + 2\sigma^2_T - \sigma^2_{T-t^\ast}}}\nonumber
\\
&&\quad{}\times\exp \biggl(-\f{A} {B + 2\sigma^2_T -
\sigma^2_{T-t^\ast}} \biggr)
\\
&&{}+ \f{N-1} {2\pi N^2 \sqrt{B + \sigma^2_T
- \sigma^2_{T-t^\ast}} \sqrt{B+\sigma^2_T
+ \sigma^2_{T-t^\ast}}}\nonumber \\
&&\quad{}\times\exp \biggl(- \f{A} {B +
\sigma^2_T + \sigma^2_{T-t^\ast}}
\biggr)\nonumber
\\
&&{}+ \f{e^{\alpha(T-t^\ast)}} {2\pi N^2 \epsilon\sqrt{B+2
\sigma^2_T}} \exp \biggl( - \f{A} {B + 2
\sigma^2_T} \biggr),\nonumber
\end{eqnarray}
where
\[
\sigma^2_s \coloneqq\f{1-e^{-2\alpha s}} {2\alpha},\qquad
A \coloneqq \bigl( e^{-\alpha T} x - y \bigr)^2,\qquad B \coloneqq
\epsilon^2 e^{-2\alpha(T-t^\ast)}.
\]
Thus, all the terms in the MSE [composed of the square
of \eqref{eq:ou-example-mean} and \eqref{eq:ou-example-var}] exhibit fairly
moderate constants, except for the last term
in \eqref{eq:ou-example-var}. Indeed, when $\alpha\gg0$, we have
$e^{\alpha(T-t^\ast)} \gg1$, unless $T - t^\ast\ll1$. In other
words, the
constant in Theorem~\ref{thr:mse-h-order} will be quite large if
$\alpha\gg
0$ and $T - t^\ast\sim1$. That observation is quite intuitive in view
of \eqref{eq:ou-example} and \eqref{eq:ou-example-reverse}: $X_{0,x}$ is
contracting to $0$ as time increases, whereas $Y_{y;T}$ is exponentially
expanding away from $y$. Thus, the probability of $X_{0,x}(t^\ast)$ and
$Y_{y;T}(T-t^\ast)$ be close to each other is very small.

%re4.22 #&#
\begin{remark}
\label{rem:ou-example-1}
Note that the last term in \eqref{eq:ou-example-var} is the term estimated
in Lemma~\ref{lem:Znm}. The constant in the lemma depends on the
constant in
Condition~\ref{ass:bound-density} for the derivatives of the transition
density $p(t,x',s,y')$ with respect to the $y'$-variable. For the
Ornstein--Uhlenbeck process, the density is given by
\[
p\bigl(t,x',s,y'\bigr) = \f{1} {2\pi
\sigma^2(s-t)} \exp \biggl(- \f{ ( e^{-\alpha
(s-t)}x'
- y' )^2} {2 \sigma^2(s-t)} \biggr).
\]
Therefore, we see that derivatives with respect to $y'$ (and, hence, the
corresponding constants) are considerably larger than derivatives with
respect to $x'$. This explains why the last term (and no other term)
in \eqref{eq:ou-example-var} causes problems for $\alpha$ large.
\end{remark}

%re4.23 #&#
\begin{remark}
\label{rem:ou-example-gen}
There is also a source of error due to the form of
$\widehat{H}_{\epsilon,N,M}$ as a fraction of two terms. The error
of an approximation
\[
\f{Q} {P} \approx\f{\widehat{Q}} {\widehat{P}} = \f{Q + \Delta Q} {P + \Delta P}
\]
of a quantity\vspace*{1pt} of interest $Q/P$ by the fraction of the approximations
$\widehat{Q}$ for $Q$ and $\widehat{P}$ for $P$ with corresponding
(absolute) errors\vadjust{\goodbreak} $\Delta Q$ and $\Delta P$ is controlled by the
\emph{relative} errors for $Q$ and $P$. Indeed, assume for simplicity that
$\Delta Q = 0$ and $Q/P = \mathcal{O}(1)$, then
\[
\f{Q} {P} - \f{\widehat{Q}} {\widehat{P}} = \mathcal{O} \biggl( \f{\Delta P /P} {1 +
\Delta P/P} \biggr),
\]
which may be close to $1$ if the relative error $\Delta P/P$ for the
denominator $P$ is large.
\end{remark}

%s5 #&#
\section{Numerical study}
\label{num}

%s5.1 #&#
\subsection{Implementation}
\label{sec:algor-impl}

Some care is necessary when implementing the forward reverse
estimators \eqref{eq:H-hat-def} and \eqref{eq:H-comp-def} for
expectations of
a functional of the diffusion bridge between two points or a point and a
subset. This especially concerns the evaluation of the double sum. Indeed,
straightforward computation would require the cost of $MN$ kernel evaluations
which would be tremendous, for example, when $M = N = 10^5$. But, fortunately,
by using kernels with an (in some sense) small support we can get
around this
difficulty as outlined below; see also \citet{MSS1} for a similar
discussion.

We here assume that the kernel $K(x)$ used in \eqref{eq:H-hat-def}
and \eqref{eq:H-comp-def}, respectively, has bounded support contained
in some
ball of radius $r$, an assumption which is easily fulfilled in
practice. For
instance, even though the Gaussian kernel $K(x)=(2\pi)^{-d/2}\exp(-|x|^2/2)$
has unbounded support, in practice $K(x)$ is negligible outside a
finite ball
(with exponential decay of the value as function of the radius).
Therefore, it
is easy to choose a ball $B_r(0)$ such that $K$ is smaller than some error
tolerance $\operatorname{const} \times\operatorname{TOL}$ outside the
ball.\footnote{Obviously, the appropriate value for $\operatorname{const}$
depends on
the size of the constants in the MSE bound.} Then, due to the small support
of $K$, the following Monte Carlo algorithm for the kernel estimator is
feasible. For simplicity, we take $N = M$.
[We present the algorithm only for the case
of \eqref{eq:H-hat-def}, the analysis being virtually equal
for \eqref{eq:H-comp-def}.] Here, the input variable $\mathcal{D}$
denotes the
grid \eqref{eq:full-grid}.
\begin{algorithm}%[!htb]
\caption{Forward-reverse algorithm}\label{alg:algorithm}
%
%{\fontsize{9.8}{11.8}{\selectfont
\begin{algorithmic}[1]
\Procedure{ForRev}{$N$, $\epsilon_N$, $a$, $\sigma$, $x$, $y$,
$\mathcal{D}$, $t^\ast$, $g$}
\State Simulate $N$ trajectories $(X^n_{s_0,x})_{n=1}^N$ of the forward
process on $s_0 < \cdots< s_K$.
\State Simulate $N$ trajectories $(Y^m_{y;T},
\mathcal{Y}_{_{y;T}}^m)_{m=1}^N$
of the reverse process on $0<\widehat{t}_1 < \cdots< \widehat{t}_L$.
\For{$m \gets1, N$}
\State Find the sub-sample
\[
\bigl\{X^{n_{k}(m)}_{s_0,x}\bigl(t^*\bigr)\dvtx k=1,
\ldots,l_m\bigr\} \coloneqq \bigl\{X^{n}_{s_0,x}
\bigl(t^*\bigr)\dvtx n=1,\ldots,N\bigr\} \cap B_{r\epsilon_N}
\bigl(Y^m_{y;T}(T)\bigr).
\]
\EndFor
\State Evaluate \eqref{eq:H-hat-def} by
{\fontsize{9}{10.8}{\selectfont
\begin{eqnarray*}\hspace*{6pt}
&&\widehat{H}_{\epsilon,M,N}\\
&&\qquad \gets\frac{\sum_{m=1}^{N} \sum_{k=1}^{l_m}
g ( X_{s_{1}}^{n_k(m)}, \ldots, X_{s_{K}}^{n_k(m)},
Y_{\widehat{t}_{L-1}}^{m}, \ldots, Y_{\widehat{t}_{1}}^{m}  )
K (
{(Y_{\widehat{t}_L}^{m} - X_{t^{\ast}}^{n_k(m)})}/{\epsilon}  )
\mathcal{Y}_{\widehat{t}_L}^{m}}{\sum_{m=1}^{N} \sum_{k=1}^{l_m} K (
{(Y_{\widehat{t}_L}^{m} - X_{t^{\ast}}^{n_k(m)})}/{\epsilon}  )
\mathcal{Y}_{\widehat{t}_L}^{m}}
\\
&&\qquad\quad{}\hspace*{2pt}\times\mathbf{1}_{{1}/{(NM)} \epsilon^{-d} \sum_{m=1}^{N} \sum_{k=1}^{l_m}
K ( {(Y_{\widehat{t}_L}^{m} - X_{t^{\ast}}^{n_k(m)})}/{\epsilon}
)
\mathcal{Y}_{\widehat{t}_L}^{m} > \overline{p}/2}.
\end{eqnarray*}}}
\EndProcedure
\end{algorithmic}
\end{algorithm}
The complexity of the simulation steps (2) and (3) in
Algorithm \ref{alg:algorithm} is $\mathcal{O} (KN d  )$ and
$\mathcal{O} (LN d  )$ elementary computations, respectively. The
size $l_m$ of the intersection in step (5) of Algorithm \ref{alg:algorithm}
is, on average, proportional to $N \epsilon_N^d \times p (s_0, x,
t^\ast,
Y^m_{y;T}(\widehat{t}_L)  )$. The search procedure itself can be done
at a cost of order $\mathcal{O}(N\log N)$ (neglecting the cost of comparison
between two integers). Thus, we get the complexity bounds summarized in
Theorem \ref{thr:complexity} below.

%th5.1 #&#
\begin{theorem}
\label{thr:complexity}
Assume that samples from the forward process $X$ and the reverse process
$(Y, \mathcal{Y})$ can be obtained at constant cost.\footnote{It is a
straightforward exercise to adjust this calculation for the case when the
corresponding stochastic differential equations need to be solved by some
numerical scheme with known rate of convergence.} Furthermore, assume that
the cost of checking for equality of integers carries negligible
cost. Then the following asymptotic bounds hold for the complexity
of Algorithm~\ref{alg:algorithm}:
\begin{itemize}
\item if $d \le4$, we choose $\epsilon_N = \mathcal{O}(N^{-1/d})$, implying
that the MSE of the output of the algorithm is $\mathcal{O}(N^{-1})$ with
a complexity $\mathcal{O}(N \log N)$;
\item if $d > 4$, we choose $\epsilon_N = \mathcal{O}(N^{-2/(4+d)})$ and
obtain an MSE of $\mathcal{O}(N^{-8/(4+d)})$ with a complexity
$\mathcal{O}(N \log N)$.
\end{itemize}
\end{theorem}

%s5.2 #&#
\subsection{Numerical examples}
\label{sec:numerical-examples}

We present two numerical studies: in the first example, the forward
process is
a two-dimensional Brownian motion, with the standard Brownian bridge as the
conditional diffusion. In the second example, we consider a Heston
model whose
stock price component is conditioned to end in a certain value. In both
examples, we actually use a Gaussian kernel
\[
K(x) = \frac{1}{(2 \pi)^{d/2}} e^{-{|x|^2}/{2}},
\]
and the simulation as well as the functional $g$ of interest are
defined on a
uniform grid $\mathcal{D} = \{ 0 = s_0 < \cdots< s_K = t^\ast= t_0 <
\cdots
< t_L = T\} $ with $s_i = i/l$ and $t_j = (K+j)/l$ for $l \in\mathbb
{N}$ and
$L+K = l$.

%ex5.2 #&#
\begin{example}
\label{ex:brownian_bridge}
We consider $X_t = B_t$, a two-dimensional standard Brownian motion, which
we condition on starting at $X_0 = 0$ and ending at $X_1 = 0$, that
is, the
conditioned diffusion is a classical two-dimensional Brownian bridge. In
particular, the reverse process $Y_t$ is also a standard Brownian motion,
and $\mathcal{Y} \equiv1$. We consider the functional
\[
g(x_1, \ldots, x_{l-1}) \coloneqq\sum
_{j=1}^2 \Biggl( \frac{1}{l-1} \sum
_{i=1}^{l-1} x^j_i
\Biggr)^2,
\]
where $x_i = (x_i^1, x_i^2) \in\mathbb{R}^2$. In this simple toy-example,
we can actually compute the true solution
\[
\mathbb{E} \bigl[ g ( X_{1/l}, \ldots, X_{(l-1)/l} )
\rrvert X_0 = X_1 = 0 \bigr] = \frac{1}{6}
\frac{l+1}{l-1}.
\]
As evaluation of the functional $g$ is cheap in this case, we use a naive
algorithm calculating the full double sum. We choose $M = N$ and
$\epsilon=
\epsilon_N = N^{-0.4}$, which still gives the rate of convergence obtained
in Theorem \ref{thr:mse-H}.
%
%f1 #&#
\begin{figure}

\includegraphics{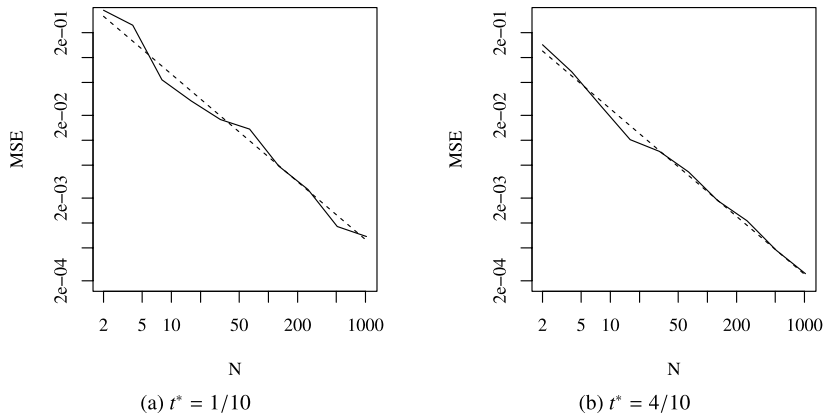}

\caption{MSE for Example \protect\ref{ex:brownian_bridge}. Dashed
lines are
reference lines proportional to $N^{-1}$.} \label{fig:ex_bb}
\end{figure}

In Figure~\ref{fig:ex_bb}, we show the results for $l=10$, with the choices
$K=1$ and $K=4$, that is, with $t^\ast= 1/10$ and $t^\ast= 4/10$,
respectively. In both case, we observe the asymptotic relation $\operatorname{MSE}
\sim N^{-1}$ predicted by Theorem \ref{thr:mse-H}. The MSE is slightly lower
when $t^\ast$ is closer to the middle of the interval $[0,1]$ [case
(b)] as
compared to the situation when $t^\ast$ is close to the boundary [case
(a)]. Intuitively, one would expect such an effect, as in the latter case
the forward process can only accumulate a considerably smaller
variance as
compared to the reverse process. However, it should be noted that the effect
is rather small.\footnote{Cf. \citet{MSS1}, where it was noted that the
variance of the forward-reverse density estimator explodes when $t^\ast
\to T$ or $t^\ast\to0$. Mathematically, this is a consequence of the
transition densities getting singular.
}
\end{example}

%ex5.3 #&#
\begin{example}
\label{ex:heston_bridge}
Let us consider the stock price $S_t$ in a Heston model: $X_t \coloneqq
(S_t, v_t)$, that is, the stock price together with its (stochastic) volatility
satisfies the stochastic differential equation
\begin{eqnarray*}
dS_t &=& \mu S_t \,dt + \sqrt{v_t}
S_t \,dB^1_t,
\\
dv_t &=& (\gamma v_t + \beta) \,dt + \xi
\sqrt{v_t} \bigl( \rho \,dB^1_t + \sqrt{1-
\rho^2} \,dB^2_t \bigr).
\end{eqnarray*}
We have
\[
a(s, x) = \pmatrix{ \mu x^1\vspace*{2pt}
\cr
\gamma x^2 +
\beta }, \qquad\sigma(s, x) = \pmatrix{ \sqrt{x^2} x^1 & 0
\vspace*{2pt}
\cr
\xi\sqrt{x^2} \rho& \xi\sqrt{x^2}
\sqrt{1-\rho^2} }.
\]
As this process is time-homogeneous, we have $\widetilde{\sigma} =
\sigma$,
and the remaining coefficients of the SDE for the reverse process are given
by
\[
\alpha(s, x) = \pmatrix{ \bigl(2 x^2 + \rho\xi- \mu \bigr)
x^1\vspace*{2pt}
\cr
(\rho\xi-\gamma) x^2 +
\xi^2 - \beta },\qquad c(s,x) = x^2 + \rho\xi- \mu- \gamma.
\]

As path-dependent functional we consider the \emph{realized variance}
of the
stock-price, that is, for a grid as above we consider
\[
g (x_1, \ldots, x_{l-1}, x_{l} ) \coloneqq\sum
_{i=1}^{l-1} \bigl( \log
\bigl(x^1_{i+1}\bigr) - \log\bigl(x^1_i
\bigr) \bigr)^2.
\]
(Dependence of the functional $g$ on the final value $y$ obviously changes
nothing in the theorems of Sections~\ref{main} and~\ref{analysis}.) We
choose $T = 1/12$ and $l = 30$. This time, however, we only condition
on the
value of the stock component at final time $T$. For the calculations,
we use
the following, arbitrarily chosen parameters: $\mu= 0.05$, $\gamma=
-0.15$, $\beta= -0.045$, $\xi= 0.3$, $\rho= -0.7$. The initial stock
price and the initial variance were set to $S_0 = 10$ and $v_0 = 0.25$,
respectively. Moreover, the realized variance was computed
conditionally on
$S_T = 12$, and we choose the standard normal density for $\varphi$,
despite Condition~\ref{ass:varphi}.

Contrary to Example \ref{ex:brownian_bridge}, we cannot produce
samples from
the exact distributions of either the forward or the reverse processes $X_t$
or $(Y_t, \mathcal{Y}_t)$. Thus, we approximate the corresponding paths
using the Euler--Maruyama scheme on a uniform grid with mesh $h = \min
(
1/360, \sqrt{0.05/N}  )$, so that the MSE for the solution of the
corresponding SDE is itself $\mathcal{O}(N^{-1})$, implying that the
asymptotic order of the MSE of our quantity of interest is not
affected by
the numerical approximation of the forward and backward processes. Moreover,
evaluation of the functional $g$ is quite costly due to the numerous calls
of the $\log$-function. Thus, we use the cut-off procedure introduced above,
so that the individual terms in the double sum are only included when the
value of the kernel $K_\epsilon$ is larger than $\eta= 10^{-3}$. The main
parameters of the forward-reverse algorithm are chosen as $M = N$ and
$\epsilon_N = (4 N)^{-0.4}$, so that we are in the regime of
Theorem \ref{thr:mse-comp}.

%f2 #&#
\begin{figure}

\includegraphics{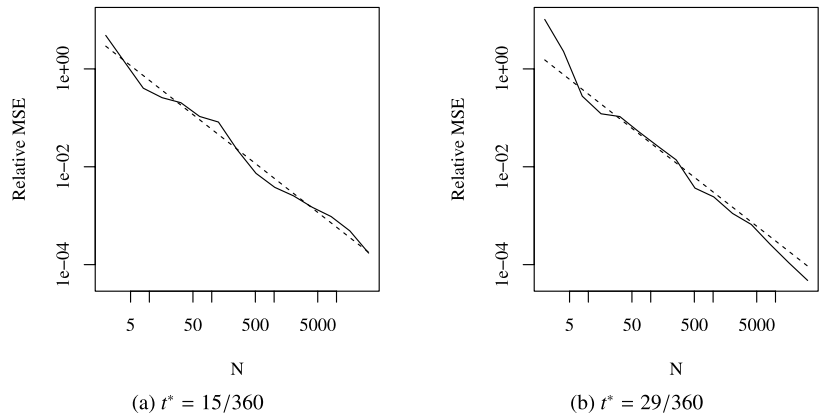}

\caption{Relative MSE for Example \protect\ref{ex:heston_bridge}.
Dashed lines are
reference lines proportional to~$N^{-1}$.}
\label{fig:ex_heston}
\end{figure}

The numerical results in Figure~\ref{fig:ex_heston} confirm the rate of
convergence for the MSE established in Theorem \ref{thr:mse-comp}. Again,
there is no significant advantage of choosing $t^\ast$ in the middle
of the
relevant interval $[0,T]$. The ``exact'' reference value was computed using
the forward-reverse algorithm with very large $N$, corresponding small
$\epsilon$ and a very fine grid for the Euler scheme. Note that
Figure~\ref{fig:ex_heston} depicts the ``relative MSE,'' that is, the MSE
normalized by the squared reference value. % Computer time was around
%five
% minutes on a Xeon Dual-Core processor.
\end{example}

% zodis "Acknowledgments" paliekamas pagal autoriu
\section*{Acknowledgments}
We are very grateful to an anonymous referee,
who has pointed out to
us the way to a much shorter and more transparent proof of the main
Theorem \ref{key}.
Moreover, the paper has profited from various
comments made by the referee,
which improved the notation and general presentation of
the paper. We are also grateful to G. N. Milstein for providing us with
enlightening references.

% imsref loaded by akundreckaite, 2014-02-20 15:08:18
%

%suskaldyti doi

\printaddresses

\end{document}